\documentclass[12pt,a4paper]{amsart}
\usepackage{PlantillaEng}
\definecolor{malaquita}{RGB}{11,218,81} 

\title{Finite generation for the group $F\left(\frac32\right)$}
\author{Jos\'e Burillo}
\address{Departament de Matem\`atiques,
Universitat Polit\`ecnica de Catalunya, Jordi Girona 1-3, 08034
Barcelona, Spain} \email{pep.burillo@upc.edu}
\thanks{The first author thanks the Spanish Ministry MICINN
through grant PID2021-126851NB-I00P for their support.}
\author{Marc Felipe}
\address{}
\email{marc.felipe.alsina@gmail.com}

\date{}

\begin{document}

\maketitle

\begin{abstract}
In this paper it is proved that the group $F\left(\frac32\right)$, a
Thompson-style group with breaks in $\ZZ\left[\frac16\right]$ but
whose slopes are restricted only to powers of $\frac32$, is finitely
generated, with a generating set of two elements.
\end{abstract}

\section*{Introduction}
Thompson's groups are a family of piecewise-linear maps of the
interval (or the circle, or the Cantor set) which have attracted a
considerable deal of attention because of their striking properties.
Some of the well-known groups of the family are the groups $F$, $T$
and $V$, see \cite{cfp}. Several families derived from these groups
have been studied along the years, such as braided versions or
higher-dimensional ones.

In \cite{BS}, Bieri and Strebel study groups of piecewise-linear maps
of intervals in great generality. Bieri and Strebel coined the
notation $G(I;A,P)$ to encode the interval of definition $I$, the
multiplicative group of slopes $P\subseteq(0,\infty)$ and the
$\ZZ[P]$-module $A\subseteq\RR$ of allowable breakpoint coordinates.
For example, Thompson's group $F$ is
$G\left([0,1];\ZZ\left[\frac12\right],\<2\>\right)$. Bieri and
Strebel groups generalize $F$ by admitting other slopes and
breakpoints. In particular, the Stein group $F(2,3)$ corresponding
to $G\left([0,1];\ZZ\left[\frac16\right],\<2,3\>\right)$ in Bieri
and Strebel notation was one of the first generalizations, studied
by Stein in \cite{stein} and also later by Wladis in \cite{wladis}.

One of the questions addressed by Bieri and Strebel in \cite{BS} is
when groups are finitely generated. They list a table of the $(A,P)$
combinations where finite generation is known. At the time of
Stein's study of $F(2,3)$, a subgroup of it, namely, the group
$F\left(\frac32\right)$, where the slopes are restricted to powers
of $\frac32$, was considered interesting, and the question arose
about whether it is finitely generated or not. The goal of this
paper is to prove that this group is indeed finitely generated, with
a set of two generators. The question of whether it is finitely presented
remains open, although in view of the formulas in \Cref{xiv2}, we suspect that possible relations
will be long and complicated.

The paper is organized as follows: after a background section, we
have four sections, defining infinite families of generators. Each
family is strictly smaller than the one before. Finally, in
\Cref{finitegen}, we prove that two generators are enough to
construct the family defined in the section right before, and hence
these two elements generate the whole group.
\section{Background}

The group $F\left(\frac32\right)=
G\left([0,1];\ZZ\left[\frac16\right],\<\frac32\>\right)$ is defined,
following \cite{BS}, as the set of orientation-preserving
homeomorphisms of $[0,1]$ that are piecewise linear, the slope of
each of the finitely many pieces is an integer power of $\frac32$,
and the points where the slope changes ---which are called
breakpoints--- have rational coordinates, with denominators having
$2$ and/or $3$ as their only prime factors. The group operation is
composition, and we write $h=f\cdot g=fg$ if $h(x)=g(f(x))$.

Elements of $F(2,3)$ are represented by tree pairs whose carets can
be binary or ternary (representing subdivisions of an interval in 2
or 3 parts), see \cite{wladis} for details. Each leaf of the tree
represents a subinterval of the type
$$\left[\frac a{2^n3^m},\frac{a+1}{2^n3^m}\right],$$
for some $a\in\ZZ$. Observe that such a leaf is located at depth
$n+m$ in the tree, because it requires $n$ binary subdivisions and
$m$ ternary ones. The particular feature of the group
$F\left(\frac32\right)$ is that since the slopes are powers of
$\frac32$, an interval of length $\frac1{2^n3^m}$ is mapped by a slope
$\left(\frac32\right)^k$ to an interval of length $\frac1{2^{n+k}3^{m-k}}$, and in
particular, the image interval is represented by a leaf \emph{at the
same depth $n+m$}. This fact characterizes the elements of
$F\left(\frac32\right)$.

\begin{rmk} The elements of $F\left(\frac32\right)$ are precisely those elements of
$F(2,3)$ such that in their tree-pair diagrams, each leaf is mapped
to another leaf at the same depth.
\end{rmk}
Observe that an obvious corollary of this fact is that every element
admits a tree-pair diagram where both trees are balanced, i.e., both
trees have all leaves in the same depth (the same depth for both
trees).
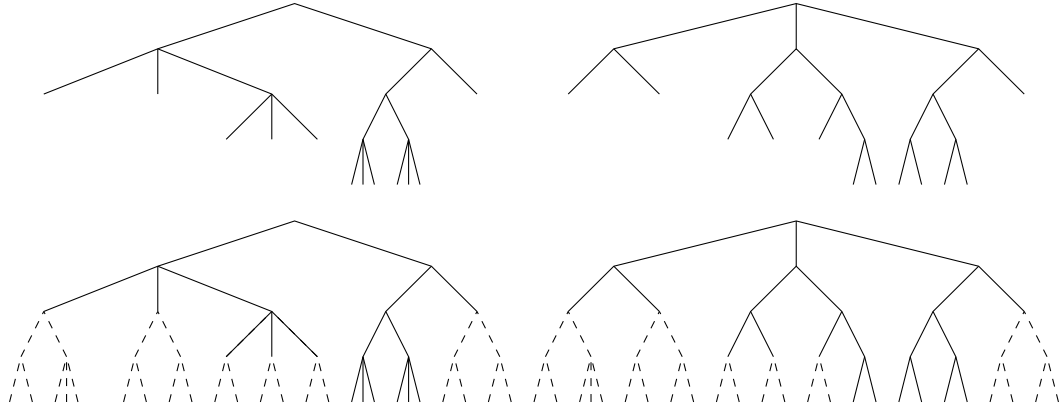
\begin{figure}[H]
    $$\begin{tikzpicture}[scale=.6]
        \node (0) at (-1,0) {};
        \node (1l) at (-4,-1) {};
        \node (1r) at (2,-1) {};
        \node (2ll) at (-6.5,-2) {};
        \node (2lm) at (-4,-2) {};
        \node (2lr) at (-1.5,-2) {};
        \node (2rl) at (1,-2) {};
        \node (2rr) at (3,-2) {};
        \node (3lrl) at (-2.5,-3) {};
        \node (3lrm) at (-1.5,-3) {};
        \node (3lrr) at (-0.5,-3) {};
        \node (3rll) at (.5,-3) {};
        \node (3rlr) at (1.5,-3) {};
        \node (4rlll) at (.25,-4) {};
        \node (4rllm) at (.5,-4) {};
        \node (4rllr) at (.75,-4) {};
        \node (4rlrl) at (1.25,-4) {};
        \node (4rlrm) at (1.5,-4) {};
        \node (4rlrr) at (1.75,-4) {};
        \draw (2ll.center)--(1l.center)--(0.center)--(1r.center)--(2rr.center);
        \draw (2lm.center)--(1l.center);
        \draw (3lrr.center)--(2lr.center)--(1l.center);
        \draw (4rlll.center)--(3rll.center)--(2rl.center)--(1r.center);
        \draw (3lrl.center)--(2lr.center);
        \draw (3lrm.center)--(2lr.center);
        \draw (4rlrr.center)--(3rlr.center)--(2rl.center);
        \draw (4rllm.center)--(3rll.center);
        \draw (4rllr.center)--(3rll.center);
        \draw (4rlrl.center)--(3rlr.center);
        \draw (4rlrm.center)--(3rlr.center);
        \node (0') at (10,0) {};
        \node (1'l) at (6,-1) {};
        \node (1'm) at (10,-1) {};
        \node (1'r) at (14,-1) {};
        \node (2'll) at (5,-2) {};
        \node (2'lr) at (7,-2) {};
        \node (2'ml) at (9,-2) {};
        \node (2'mr) at (11,-2) {};
        \node (2'rl) at (13,-2) {};
        \node (2'rr) at (15,-2) {};
        \node (3'mll) at (8.5,-3) {};
        \node (3'mlr) at (9.5,-3) {};
        \node (3'mrl) at (10.5,-3) {};
        \node (3'mrr) at (11.5,-3) {};
        \node (3'rll) at (12.5,-3) {};
        \node (3'rlr) at (13.5,-3) {};
        \node (4'mrrl) at (11.25,-4) {};
        \node (4'mrrr) at (11.75,-4) {};
        \node (4'rlll) at (12.25,-4) {};
        \node (4'rllr) at (12.75,-4) {};
        \node (4'rlrl) at (13.25,-4) {};
        \node (4'rlrr) at (13.75,-4) {};
        \draw (2'll.center)--(1'l.center)--(0'.center)--(1'r.center)--(2'rr.center);
        \draw (1'm.center)--(0'.center);
        \draw (2'lr.center)--(1'l.center);
        \draw (3'mll.center)--(2'ml.center)--(1'm.center)--(2'mr.center)--(3'mrr.center)--(4'mrrr.center);
        \draw (4'rlll.center)--(3'rll.center)--(2'rl.center)--(1'r.center);
        \draw (3'mlr.center)--(2'ml.center);
        \draw (3'mrl.center)--(2'mr.center);
        \draw (4'rlrr.center)--(3'rlr.center)--(2'rl.center);
        \draw (4'mrrl.center)--(3'mrr.center);
        \draw (4'rllr.center)--(3'rll.center);
        \draw (4'rlrl.center)--(3'rlr.center);
    \end{tikzpicture}$$
    $$\begin{tikzpicture}[scale=.6]
        \node (0) at (-1,0) {};
        \node (1l) at (-4,-1) {};
        \node (1r) at (2,-1) {};
        \node (2ll) at (-6.5,-2) {};
        \node (2lm) at (-4,-2) {};
        \node (2lr) at (-1.5,-2) {};
        \node (2rl) at (1,-2) {};
        \node (2rr) at (3,-2) {};
        \node (3lll) at (-7,-3) {};
        \node (3llr) at (-6,-3) {};
        \node (3lml) at (-4.5,-3) {};
        \node (3lmr) at (-3.5,-3) {};
        \node (3lrl) at (-2.5,-3) {};
        \node (3lrm) at (-1.5,-3) {};
        \node (3lrr) at (-0.5,-3) {};
        \node (3rll) at (.5,-3) {};
        \node (3rlr) at (1.5,-3) {};
        \node (3rrl) at (2.5,-3) {};
        \node (3rrr) at (3.5,-3) {};
        \node (4llll) at (-7.25,-4) {};
        \node (4lllr) at (-6.75,-4) {};
        \node (4llrl) at (-6.25,-4) {};
        \node (4llrm) at (-6,-4) {};
        \node (4llrr) at (-5.75,-4) {};
        \node (4lmll) at (-4.75,-4) {};
        \node (4lmlr) at (-4.25,-4) {};
        \node (4lmrl) at (-3.75,-4) {};
        \node (4lmrr) at (-3.25,-4) {};
        \node (4lrll) at (-2.75,-4) {};
        \node (4lrlr) at (-2.25,-4) {};
        \node (4lrml) at (-1.75,-4) {};
        \node (4lrmr) at (-1.25,-4) {};
        \node (4lrrl) at (-.75,-4) {};
        \node (4lrrr) at (-.25,-4) {};
        \node (4rlll) at (.25,-4) {};
        \node (4rllm) at (.5,-4) {};
        \node (4rllr) at (.75,-4) {};
        \node (4rlrl) at (1.25,-4) {};
        \node (4rlrm) at (1.5,-4) {};
        \node (4rlrr) at (1.75,-4) {};
        \node (4rrll) at (2.25,-4) {};
        \node (4rrlr) at (2.75,-4) {};
        \node (4rrrl) at (3.25,-4) {};
        \node (4rrrr) at (3.75,-4) {};
        \draw (2ll.center)--(1l.center)--(0.center)--(1r.center)--(2rr.center);
        \draw (2lm.center)--(1l.center);
        \draw (3lrr.center)--(2lr.center)--(1l.center);
        \draw (4rlll.center)--(3rll.center)--(2rl.center)--(1r.center);
        \draw (3lrl.center)--(2lr.center);
        \draw (3lrm.center)--(2lr.center);
        \draw (4rlrr.center)--(3rlr.center)--(2rl.center);
        \draw (4rllm.center)--(3rll.center);
        \draw (4rllr.center)--(3rll.center);
        \draw (4rlrl.center)--(3rlr.center);
        \draw (4rlrm.center)--(3rlr.center);
        \draw[dashed] (3lll.center)--(2ll.center);
        \draw[dashed] (3llr.center)--(2ll.center);
        \draw[dashed] (3lml.center)--(2lm.center);
        \draw[dashed] (3lmr.center)--(2lm.center);
        \draw[dashed] (3lrl.center)--(2lr.center);
        \draw[dashed] (3lrr.center)--(2lr.center);
        \draw[dashed] (3rrl.center)--(2rr.center);
        \draw[dashed] (3rrr.center)--(2rr.center);
        \draw[dashed] (4llll.center)--(3lll.center);
        \draw[dashed] (4lllr.center)--(3lll.center);
        \draw[dashed] (4llrl.center)--(3llr.center);
        \draw[dashed] (4llrm.center)--(3llr.center);
        \draw[dashed] (4llrr.center)--(3llr.center);
        \draw[dashed] (4lmll.center)--(3lml.center);
        \draw[dashed] (4lmlr.center)--(3lml.center);
        \draw[dashed] (4lmrl.center)--(3lmr.center);
        \draw[dashed] (4lmrr.center)--(3lmr.center);
        \draw[dashed] (4lrll.center)--(3lrl.center);
        \draw[dashed] (4lrlr.center)--(3lrl.center);
        \draw[dashed] (4lrml.center)--(3lrm.center);
        \draw[dashed] (4lrmr.center)--(3lrm.center);
        \draw[dashed] (4lrrl.center)--(3lrr.center);
        \draw[dashed] (4lrrr.center)--(3lrr.center);
        \draw[dashed] (4rrll.center)--(3rrl.center);
        \draw[dashed] (4rrlr.center)--(3rrl.center);
        \draw[dashed] (4rrrl.center)--(3rrr.center);
        \draw[dashed] (4rrrr.center)--(3rrr.center);
        \node (0') at (10,0) {};
        \node (1'l) at (6,-1) {};
        \node (1'm) at (10,-1) {};
        \node (1'r) at (14,-1) {};
        \node (2'll) at (5,-2) {};
        \node (2'lr) at (7,-2) {};
        \node (2'ml) at (9,-2) {};
        \node (2'mr) at (11,-2) {};
        \node (2'rl) at (13,-2) {};
        \node (2'rr) at (15,-2) {};
        \node (3'lll) at (4.5,-3) {};
        \node (3'llr) at (5.5,-3) {};
        \node (3'lrl) at (6.5,-3) {};
        \node (3'lrr) at (7.5,-3) {};
        \node (3'mll) at (8.5,-3) {};
        \node (3'mlr) at (9.5,-3) {};
        \node (3'mrl) at (10.5,-3) {};
        \node (3'mrr) at (11.5,-3) {};
        \node (3'rll) at (12.5,-3) {};
        \node (3'rlr) at (13.5,-3) {};
        \node (3'rrl) at (14.5,-3) {};
        \node (3'rrr) at (15.5,-3) {};
        \node (4'llll) at (4.25,-4) {};
        \node (4'lllr) at (4.75,-4) {};
        \node (4'llrl) at (5.25,-4) {};
        \node (4'llrm) at (5.5,-4) {};
        \node (4'llrr) at (5.75,-4) {};
        \node (4'lrll) at (6.25,-4) {};
        \node (4'lrlr) at (6.75,-4) {};
        \node (4'lrrl) at (7.25,-4) {};
        \node (4'lrrr) at (7.75,-4) {};
        \node (4'mlll) at (8.25,-4) {};
        \node (4'mllr) at (8.75,-4) {};
        \node (4'mlrl) at (9.25,-4) {};
        \node (4'mlrr) at (9.75,-4) {};
        \node (4'mrll) at (10.25,-4) {};
        \node (4'mrlr) at (10.75,-4) {};
        \node (4'mrrl) at (11.25,-4) {};
        \node (4'mrrr) at (11.75,-4) {};
        \node (4'rlll) at (12.25,-4) {};
        \node (4'rllr) at (12.75,-4) {};
        \node (4'rlrl) at (13.25,-4) {};
        \node (4'rlrr) at (13.75,-4) {};
        \node (4'rrll) at (14.25,-4) {};
        \node (4'rrlr) at (14.75,-4) {};
        \node (4'rrrl) at (15.25,-4) {};
        \node (4'rrrr) at (15.75,-4) {};
        \draw (2'll.center)--(1'l.center)--(0'.center)--(1'r.center)--(2'rr.center);
        \draw (1'm.center)--(0'.center);
        \draw (2'lr.center)--(1'l.center);
        \draw (3'mll.center)--(2'ml.center)--(1'm.center)--(2'mr.center)--(3'mrr.center)--(4'mrrr.center);
        \draw (4'rlll.center)--(3'rll.center)--(2'rl.center)--(1'r.center);
        \draw (3'mlr.center)--(2'ml.center);
        \draw (3'mrl.center)--(2'mr.center);
        \draw (4'rlrr.center)--(3'rlr.center)--(2'rl.center);
        \draw (4'mrrl.center)--(3'mrr.center);
        \draw (4'rllr.center)--(3'rll.center);
        \draw (4'rlrl.center)--(3'rlr.center);
        \draw[dashed] (3'lll.center)--(2'll.center);
        \draw[dashed] (3'llr.center)--(2'll.center);
        \draw[dashed] (3'lrl.center)--(2'lr.center);
        \draw[dashed] (3'lrr.center)--(2'lr.center);
        \draw[dashed] (3'rrl.center)--(2'rr.center);
        \draw[dashed] (3'rrr.center)--(2'rr.center);
        \draw[dashed] (4'llll.center)--(3'lll.center);
        \draw[dashed] (4'lllr.center)--(3'lll.center);
        \draw[dashed] (4'llrl.center)--(3'llr.center);
        \draw[dashed] (4'llrm.center)--(3'llr.center);
        \draw[dashed] (4'llrr.center)--(3'llr.center);
        \draw[dashed] (4'lrll.center)--(3'lrl.center);
        \draw[dashed] (4'lrlr.center)--(3'lrl.center);
        \draw[dashed] (4'lrrl.center)--(3'lrr.center);
        \draw[dashed] (4'lrrr.center)--(3'lrr.center);
        \draw[dashed] (4'mlll.center)--(3'mll.center);
        \draw[dashed] (4'mllr.center)--(3'mll.center);
        \draw[dashed] (4'mlrl.center)--(3'mlr.center);
        \draw[dashed] (4'mlrr.center)--(3'mlr.center);
        \draw[dashed] (4'mrll.center)--(3'mrl.center);
        \draw[dashed] (4'mrlr.center)--(3'mrl.center);
        \draw[dashed] (4'rrll.center)--(3'rrl.center);
        \draw[dashed] (4'rrlr.center)--(3'rrl.center);
        \draw[dashed] (4'rrrl.center)--(3'rrr.center);
        \draw[dashed] (4'rrrr.center)--(3'rrr.center);
    \end{tikzpicture}$$
    \caption{Above, we see a tree-pair diagram of an element of $F\left(\frac32\right)$. From left to right, the leaves are at depths $2, 2, 3, 3, 3, 4, 4, 4, 4, 4, 4, 2$ on both trees. Below, some carets have been added to make a balanced tree-pair diagram of the same element, whose trees have 25 leaves each, all at depth 4. Note that both 2-carets and 3-carets can be used arbitrarily to construct the balanced tree, as long as the same arrangement of carets is used for the other tree.}
\end{figure}
Two of the easiest elements of $F\left(\frac32\right)$, having only
two rows of carets, are the elements $l$ and $r$, given by
$$
l(t)=\left\{\begin{array}{ll} \frac{3t}2&\text{ if
}t\in[0,\frac29]\\
\vspace{-3mm}\\
t+\frac19&\text{ if }t\in[\frac29,\frac13]\\
\vspace{-3mm}\\
\frac{2t}3+\frac29&\text{ if }t\in[\frac13,\frac23]\\
\vspace{-3mm}\\
t&\text{ if }t\in[\frac23,1]\end{array}\right. \qquad
r(t)=\left\{\begin{array}{ll} t&\text{ if
}t\in[0,\frac13]\\
\vspace{-3mm}\\
\frac {3t}2-\frac16&\text{ if }t\in[\frac13,\frac59]\\
\vspace{-3mm}\\
t+\frac19&\text{ if }t\in[\frac59,\frac23]\\
\vspace{-3mm}\\
\frac {2t}3+\frac13&\text{ if }t\in[\frac23,1]\end{array}\right.
$$

\begin{figure}[H]
        \begin{align*}\begin{tikzpicture}[scale=5]
            \node (l) at (-.15,1.15) {$l\text{:}$};
            \draw[thick] (0,0) -- (2/9,1/3) -- (1/3,4/9) -- (2/3,2/3) -- (1,1);
            \draw[ultra thick] (-.05,0) -- (1.05,0);
            \draw[ultra thick] (0,-.05) -- (0,1.05);
            \draw[dashed,color=gray] (1,0) -- (1,1.05);
            \node[anchor=north] (x1) at (1,0) {$1$};
            \draw[dashed,color=gray] (0,1) -- (1.05,1);
            \node[anchor=east] (y1) at (0,1) {$1$};
            \draw[dashed,color=gray] (2/3,0) -- (2/3,1.05);
            \node[anchor=north] (x2/3) at (2/3,0) {$\frac23$};
            \draw[dashed,color=gray] (0,2/3) -- (1.05,2/3);
            \node[anchor=east] (y2/3) at (0,2/3) {$\frac23$};
            \draw[dashed,color=gray] (1/3,0) -- (1/3,1.05);
            \node[anchor=north] (x1/3) at (1/3,0) {$\frac13$};
            \draw[dashed,color=gray] (0,1/3) -- (1.05,1/3);
            \node[anchor=east] (y1/3) at (0,1/3) {$\frac13$};
            \draw[dotted] (1/9,0) -- (1/9,1.05);
            \draw[dotted] (2/9,0) -- (2/9,1.05);
            \draw[dotted] (4/9,0) -- (4/9,1.05);
            \draw[dotted] (5/9,0) -- (5/9,1.05);
            \draw[dotted] (7/9,0) -- (7/9,1.05);
            \draw[dotted] (8/9,0) -- (8/9,1.05);
            \draw[dotted] (0,1/9) -- (1.05,1/9);
            \draw[dotted] (0,2/9) -- (1.05,2/9);
            \draw[dotted] (0,4/9) -- (1.05,4/9);
            \draw[dotted] (0,5/9) -- (1.05,5/9);
            \draw[dotted] (0,7/9) -- (1.05,7/9);
            \draw[dotted] (0,8/9) -- (1.05,8/9);
        \end{tikzpicture}&\qquad\begin{tikzpicture}[scale=5]
            \node (r) at (-.15,1.15) {$r\text{:}$};
            \draw[thick] (0,0) -- (1/3,1/3) -- (5/9,2/3) -- (2/3,7/9) -- (1,1);
            \draw[ultra thick] (-.05,0) -- (1.05,0);
            \draw[ultra thick] (0,-.05) -- (0,1.05);
            \draw[dashed,color=gray] (1,0) -- (1,1.05);
            \node[anchor=north] (x1) at (1,0) {$1$};
            \draw[dashed,color=gray] (0,1) -- (1.05,1);
            \node[anchor=east] (y1) at (0,1) {$1$};
            \draw[dashed,color=gray] (2/3,0) -- (2/3,1.05);
            \node[anchor=north] (x2/3) at (2/3,0) {$\frac23$};
            \draw[dashed,color=gray] (0,2/3) -- (1.05,2/3);
            \node[anchor=east] (y2/3) at (0,2/3) {$\frac23$};
            \draw[dashed,color=gray] (1/3,0) -- (1/3,1.05);
            \node[anchor=north] (x1/3) at (1/3,0) {$\frac13$};
            \draw[dashed,color=gray] (0,1/3) -- (1.05,1/3);
            \node[anchor=east] (y1/3) at (0,1/3) {$\frac13$};
            \draw[dotted] (1/9,0) -- (1/9,1.05);
            \draw[dotted] (2/9,0) -- (2/9,1.05);
            \draw[dotted] (4/9,0) -- (4/9,1.05);
            \draw[dotted] (5/9,0) -- (5/9,1.05);
            \draw[dotted] (7/9,0) -- (7/9,1.05);
            \draw[dotted] (8/9,0) -- (8/9,1.05);
            \draw[dotted] (0,1/9) -- (1.05,1/9);
            \draw[dotted] (0,2/9) -- (1.05,2/9);
            \draw[dotted] (0,4/9) -- (1.05,4/9);
            \draw[dotted] (0,5/9) -- (1.05,5/9);
            \draw[dotted] (0,7/9) -- (1.05,7/9);
            \draw[dotted] (0,8/9) -- (1.05,8/9);
        \end{tikzpicture}\\
        \raisebox{-\height/2}{\begin{tikzpicture}
        \node (l) at (-1,0) {$l\text{:}$};
        \node (0) at (0,0) {};
        \node (1l) at (-1,-1) {};
        \node (1m) at (0,-1) {};
        \node (1r) at (1,-1) {};
        \node (2l) at (-1.25,-2) {};
        \node (2m) at (-1,-2) {};
        \node (2r) at (-.75,-2) {};
        \node (2L) at (-.25,-2) {};
        \node (2R) at (.25,-2) {};
        \draw (2l.center)--(1l.center)--(0.center)--(1r.center);
        \draw (1m.center)--(0.center);
        \draw (2r.center)--(1l.center);
        \draw (2m.center)--(1l.center);
        \draw (2L.center)--(1m.center);
        \draw (2R.center)--(1m.center);
        \node (0') at (3,0) {};
        \node (1l') at (2,-1) {};
        \node (1m') at (3,-1) {};
        \node (1r') at (4,-1) {};
        \node (2l') at (1.75,-2) {};
        \node (2r') at (2.25,-2) {};
        \node (2L') at (2.75,-2) {};
        \node (2M') at (3,-2) {};
        \node (2R') at (3.25,-2) {};
        \draw (2l'.center)--(1l'.center)--(0'.center)--(1r'.center);
        \draw (1m'.center)--(0'.center);
        \draw (2r'.center)--(1l'.center);
        \draw (2L'.center)--(1m'.center);
        \draw (2M'.center)--(1m'.center);
        \draw (2R'.center)--(1m'.center);
    \end{tikzpicture}}&\qquad\raisebox{-\height/2}{\begin{tikzpicture}
        \node (r) at (-1,0) {$r\text{:}$};
        \node (0) at (0,0) {};
        \node (1l) at (-1,-1) {};
        \node (1m) at (0,-1) {};
        \node (1r) at (1,-1) {};
        \node (2l) at (-.25,-2) {};
        \node (2m) at (0,-2) {};
        \node (2r) at (.25,-2) {};
        \node (2R) at (1.25,-2) {};
        \node (2L) at (.75,-2) {};
        \draw (1l.center)--(0.center)--(1r.center)--(2R.center);
        \draw (1m.center)--(0.center);
        \draw (2l.center)--(1m.center);
        \draw (2m.center)--(1m.center);
        \draw (2r.center)--(1m.center);
        \draw (2L.center)--(1r.center);
        \node (0') at (3,0) {};
        \node (1l') at (2,-1) {};
        \node (1m') at (3,-1) {};
        \node (1r') at (4,-1) {};
        \node (2l') at (2.75,-2) {};
        \node (2r') at (3.25,-2) {};
        \node (2L') at (3.75,-2) {};
        \node (2M') at (4,-2) {};
        \node (2R') at (4.25,-2) {};
        \draw (1l'.center)--(0'.center)--(1r'.center)--(2R'.center);
        \draw (1m'.center)--(0'.center);
        \draw (2l'.center)--(1m'.center);
        \draw (2r'.center)--(1m'.center);
        \draw (2L'.center)--(1r'.center);
        \draw (2M'.center)--(1r'.center);
    \end{tikzpicture}}\end{align*}
        \caption{The generators $l$ and $r$ and their tree-pair diagram representation.}
        \label{lir}
    \end{figure}
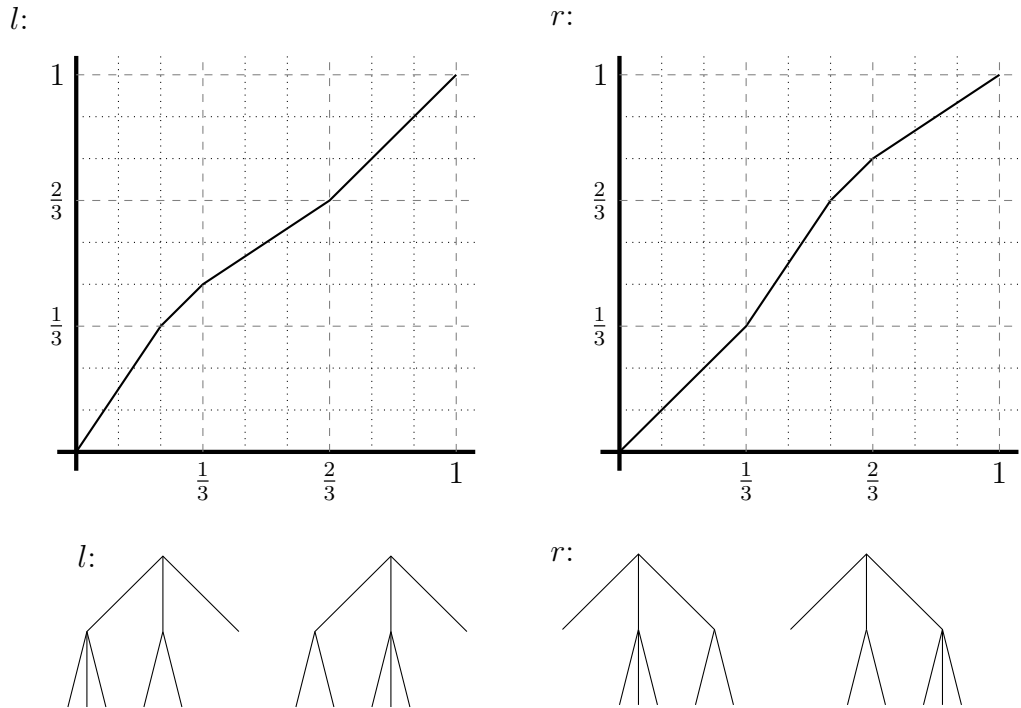

See \Cref{lir} for the graphs and the tree-pair diagrams for $l$ and
$r$.

The main theorem of this paper is then the following.

\begin{thm}\label{fg}
The group $F\left(\frac32\right)$ is generated by the maps $l$ and $r$.
\end{thm}

The proof of this theorem will take the rest of the paper. To this
end, we will give several families of generators, which we call
the \textit{humanoids}, the \textit{serpents}, the
\textit{centipedes} and \textit{stick bugs} for reasons that will
become apparent later. We will prove that each family is a
generating subset of the previous one and show at the end that $l$
and $r$ will be enough to generate the last family.

\section{Humanoids}

    We start with the first of our plethora of generating families, the humanoids.

    \begin{dfn}
    An element is called \emph{humanoid} if it (or its inverse) admits a tree-pair diagram of the form shown in \Cref{humanoid}, where $T$ is some tree and the black dots represent consecutive leaves of $T$ that are at the same depth.\\
    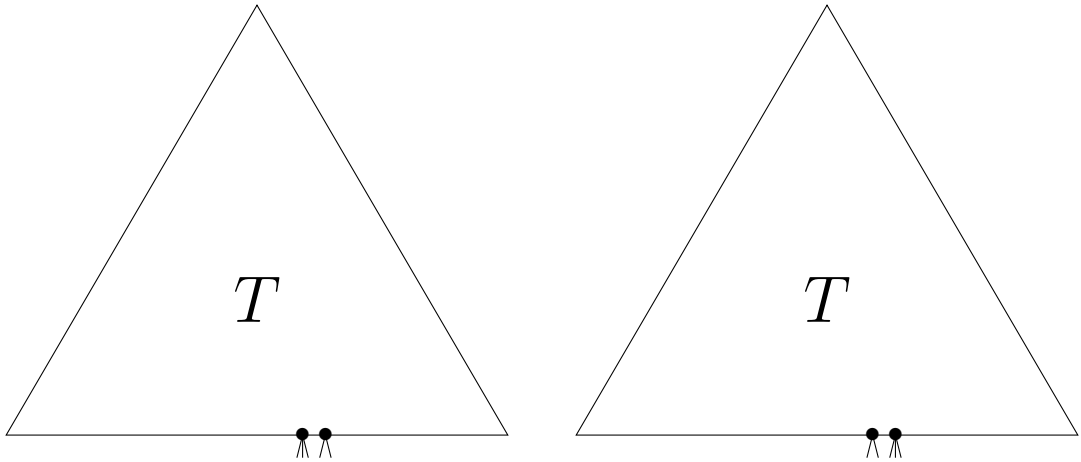
\begin{figure}[H]
    $$\begin{tikzpicture}[scale=.3]
        \node (t1) at (0,19) {};
        \node (l1) at (-11,0) {};
        \node (r1) at (11,0) {};
        \node (c1) at (0,6) {\Huge $T$};
        \draw (r1.center)--(t1.center)--(l1.center)--cycle;
        \node (a1) at (2,0) {$\bullet$};
        \node (b1) at (3,0) {$\bullet$};
        \node (d1) at (2.25,-1) {};
        \node (e1) at (2,-1) {};
        \node (f1) at (1.75,-1) {};
        \node (g1) at (3.25,-1) {};
        \node (h1) at (2.75,-1) {};
        \draw (a1.center)--(d1.center);
        \draw (a1.center)--(e1.center);
        \draw (a1.center)--(f1.center);
        \draw (b1.center)--(g1.center);
        \draw (b1.center)--(h1.center);
        \node (t2) at (25,19) {};
        \node (l2) at (14,0) {};
        \node (r2) at (36,0) {};
        \node (c2) at (25,6) {\Huge $T$};
        \draw (r2.center)--(t2.center)--(l2.center)--cycle;
        \node (a2) at (27,0) {$\bullet$};
        \node (b2) at (28,0) {$\bullet$};
        \node (d2) at (27.25,-1) {};
        \node (e2) at (26.75,-1) {};
        \node (f2) at (28.25,-1) {};
        \node (g2) at (28,-1) {};
        \node (h2) at (27.75,-1) {};
        \draw (a2.center)--(d2.center);
        \draw (a2.center)--(e2.center);
        \draw (b2.center)--(f2.center);
        \draw (b2.center)--(g2.center);
        \draw (b2.center)--(h2.center);
    \end{tikzpicture}$$
    \caption{A humanoid tree-pair diagram.}
    \label{humanoid}
    \end{figure}
    \end{dfn}

    \begin{prp}\label{humgen}
    Humanoids generate $F\left(\frac32\right)$.
    \end{prp}

    In order to prove \Cref{humgen}, we will make use of the following lemma.

\begin{figure}[H]
    $$\left(\raisebox{-\height/2}{\begin{tikzpicture} [scale=.7]
        \node (0) at (0,0) {};
        \node (1l) at (-2,-1) {};
        \node (1m) at (0,-1) {};
        \node (1r) at (2,-1) {};
        \node (2lr) at (-1.5,-2) {};
        \node (2lm) at (-2,-2) {};
        \node (2ll) at (-2.5,-2) {};
        \node (2mr) at (.5,-2) {};
        \node (2ml) at (-.5,-2) {};
        \node (2rr) at (2.5,-2) {};
        \node (2rl) at (1.5,-2) {};
        \node (3lrl) at (-1.6,-3) {};
        \node (3lrm) at (-1.5,-3) {};
        \node (3lrr) at (-1.4,-3) {};
        \node (3mll) at (-.6,-3) {};
        \node (3mlr) at (-.4,-3) {};
        \draw (2ll.center)--(1l.center)--(0.center)--(1r.center)--(2rr.center);
        \draw (1m.center)--(0.center);
        \draw (2lm.center)--(1l.center);
        \draw (2lr.center)--(1l.center);
        \draw (2ml.center)--(1m.center);
        \draw (2mr.center)--(1m.center);
        \draw (2rl.center)--(1r.center);
        \draw (3mll.center)--(2ml.center);
        \draw (3mlr.center)--(2ml.center);
        \draw (3lrl.center)--(2lr.center);
        \draw (3lrm.center)--(2lr.center);
        \draw (3lrr.center)--(2lr.center);
    \end{tikzpicture}}\raisebox{-\height/2}{\begin{tikzpicture} [scale=.7]
        \node (0) at (0,0) {};
        \node (1l) at (-2,-1) {};
        \node (1m) at (0,-1) {};
        \node (1r) at (2,-1) {};
        \node (2lr) at (-1.5,-2) {};
        \node (2lm) at (-2,-2) {};
        \node (2ll) at (-2.5,-2) {};
        \node (2mr) at (.5,-2) {};
        \node (2ml) at (-.5,-2) {};
        \node (2rr) at (2.5,-2) {};
        \node (2rl) at (1.5,-2) {};
        \node (3lrl) at (-1.6,-3) {};
        \node (3lrr) at (-1.4,-3) {};
        \node (3mll) at (-.6,-3) {};
        \node (3mlm) at (-.5,-3) {};
        \node (3mlr) at (-.4,-3) {};
        \draw (2ll.center)--(1l.center)--(0.center)--(1r.center)--(2rr.center);
        \draw (1m.center)--(0.center);
        \draw (2lm.center)--(1l.center);
        \draw (2lr.center)--(1l.center);
        \draw (2ml.center)--(1m.center);
        \draw (2mr.center)--(1m.center);
        \draw (2rl.center)--(1r.center);
        \draw (3mll.center)--(2ml.center);
        \draw (3mlm.center)--(2ml.center);
        \draw (3mlr.center)--(2ml.center);
        \draw (3lrl.center)--(2lr.center);
        \draw (3lrr.center)--(2lr.center);
        \draw[dashed] (-2.8,0.3)--(-2.8,-2.3)--(-1.75,-2.3)--(-1.75,-3.3)--(-.25,-3.3)--(-.25,-2.3)--(2.8,-2.3)--(2.8,0.3)--cycle;
    \end{tikzpicture}}\right)\left(\raisebox{-\height/2}{\begin{tikzpicture} [scale=.7]
        \node (0) at (0,0) {};
        \node (1l) at (-2,-1) {};
        \node (1m) at (0,-1) {};
        \node (1r) at (2,-1) {};
        \node (2lr) at (-1.5,-2) {};
        \node (2lm) at (-2,-2) {};
        \node (2ll) at (-2.5,-2) {};
        \node (2mr) at (.5,-2) {};
        \node (2ml) at (-.5,-2) {};
        \node (2rr) at (2.5,-2) {};
        \node (2rl) at (1.5,-2) {};
        \node (3lll) at (-2.6,-3) {};
        \node (3llr) at (-2.4,-3) {};
        \node (3lml) at (-2.1,-3) {};
        \node (3lmm) at (-2,-3) {};
        \node (3lmr) at (-1.9,-3) {};
        \node (3lrl) at (-1.6,-3) {};
        \node (3lrr) at (-1.4,-3) {};
        \node (3mll) at (-.6,-3) {};
        \node (3mlm) at (-.5,-3) {};
        \node (3mlr) at (-.4,-3) {};
        \node (3mrl) at (.4,-3) {};
        \node (3mrr) at (.6,-3) {};
        \node (3rll) at (1.4,-3) {};
        \node (3rlr) at (1.6,-3) {};
        \node (3rrl) at (2.4,-3) {};
        \node (3rrr) at (2.6,-3) {};
        \draw (2ll.center)--(1l.center)--(0.center)--(1r.center)--(2rr.center);
        \draw (1m.center)--(0.center);
        \draw (2lm.center)--(1l.center);
        \draw (2lr.center)--(1l.center);
        \draw (2ml.center)--(1m.center);
        \draw (2mr.center)--(1m.center);
        \draw (2rl.center)--(1r.center);
        \draw (3lll.center)--(2ll.center);
        \draw (3llr.center)--(2ll.center);
        \draw (3lml.center)--(2lm.center);
        \draw (3lmm.center)--(2lm.center);
        \draw (3lmr.center)--(2lm.center);
        \draw (3mll.center)--(2ml.center);
        \draw (3mlm.center)--(2ml.center);
        \draw (3mlr.center)--(2ml.center);
        \draw (3mrl.center)--(2mr.center);
        \draw (3mrr.center)--(2mr.center);
        \draw (3lrl.center)--(2lr.center);
        \draw (3lrr.center)--(2lr.center);
        \draw (3rll.center)--(2rl.center);
        \draw (3rlr.center)--(2rl.center);
        \draw (3rrl.center)--(2rr.center);
        \draw (3rrr.center)--(2rr.center);
        \draw[dashed] (-2.8,0.3)--(-2.8,-2.3)--(-1.75,-2.3)--(-1.75,-3.3)--(-.25,-3.3)--(-.25,-2.3)--(2.8,-2.3)--(2.8,0.3)--cycle;
    \end{tikzpicture}}\raisebox{-\height/2}{\begin{tikzpicture} [scale=.7]
        \node (0) at (0,0) {};
        \node (1l) at (-1,-1) {};
        \node (1r) at (1,-1) {};
        \node (2ll) at (-1.5,-2) {};
        \node (2lm) at (-1,-2) {};
        \node (2lr) at (-.5,-2) {};
        \node (2rl) at (.5,-2) {};
        \node (2rm) at (1,-2) {};
        \node (2rr) at (1.5,-2) {};
        \node (3lll) at (-1.6,-3) {};
        \node (3llm) at (-1.5,-3) {};
        \node (3llr) at (-1.4,-3) {};
        \node (3lml) at (-1.1,-3) {};
        \node (3lmr) at (-.9,-3) {};
        \node (3lrl) at (-.6,-3) {};
        \node (3lrm) at (-.5,-3) {};
        \node (3lrr) at (-.4,-3) {};
        \node (3rrr) at (1.6,-3) {};
        \node (3rrm) at (1.5,-3) {};
        \node (3rrl) at (1.4,-3) {};
        \node (3rmr) at (1.1,-3) {};
        \node (3rml) at (.9,-3) {};
        \node (3rlr) at (.6,-3) {};
        \node (3rlm) at (.5,-3) {};
        \node (3rll) at (.4,-3) {};
        \draw (2ll.center)--(1l.center)--(0.center)--(1r.center)--(2rr.center);
        \draw (2lm.center)--(1l.center);
        \draw (2lr.center)--(1l.center);
        \draw (2rl.center)--(1r.center);
        \draw (2rm.center)--(1r.center);
        \draw (3lll.center)--(2ll.center);
        \draw (3llm.center)--(2ll.center);
        \draw (3llr.center)--(2ll.center);
        \draw (3lml.center)--(2lm.center);
        \draw (3lmr.center)--(2lm.center);
        \draw (3lrl.center)--(2lr.center);
        \draw (3lrm.center)--(2lr.center);
        \draw (3lrr.center)--(2lr.center);
        \draw (3rrr.center)--(2rr.center);
        \draw (3rrm.center)--(2rr.center);
        \draw (3rrl.center)--(2rr.center);
        \draw (3rmr.center)--(2rm.center);
        \draw (3rml.center)--(2rm.center);
        \draw (3rlr.center)--(2rl.center);
        \draw (3rlm.center)--(2rl.center);
        \draw (3rll.center)--(2rl.center);
    \end{tikzpicture}}\right)=$$
    $$=\raisebox{-\height/2}{\begin{tikzpicture} [scale=.7]
        \node (0) at (0,0) {};
        \node (1l) at (-2,-1) {};
        \node (1m) at (0,-1) {};
        \node (1r) at (2,-1) {};
        \node (2lr) at (-1.5,-2) {};
        \node (2lm) at (-2,-2) {};
        \node (2ll) at (-2.5,-2) {};
        \node (2mr) at (.5,-2) {};
        \node (2ml) at (-.5,-2) {};
        \node (2rr) at (2.5,-2) {};
        \node (2rl) at (1.5,-2) {};
        \node (3lll) at (-2.6,-3) {};
        \node (3llr) at (-2.4,-3) {};
        \node (3lml) at (-2.1,-3) {};
        \node (3lmm) at (-2,-3) {};
        \node (3lmr) at (-1.9,-3) {};
        \node (3lrl) at (-1.6,-3) {};
        \node (3lrm) at (-1.5,-3) {};
        \node (3lrr) at (-1.4,-3) {};
        \node (3mll) at (-.6,-3) {};
        \node (3mlr) at (-.4,-3) {};
        \node (3mrl) at (.4,-3) {};
        \node (3mrr) at (.6,-3) {};
        \node (3rll) at (1.4,-3) {};
        \node (3rlr) at (1.6,-3) {};
        \node (3rrl) at (2.4,-3) {};
        \node (3rrr) at (2.6,-3) {};
        \draw (2ll.center)--(1l.center)--(0.center)--(1r.center)--(2rr.center);
        \draw (1m.center)--(0.center);
        \draw (2lm.center)--(1l.center);
        \draw (2lr.center)--(1l.center);
        \draw (2ml.center)--(1m.center);
        \draw (2mr.center)--(1m.center);
        \draw (2rl.center)--(1r.center);
        \draw (3lll.center)--(2ll.center);
        \draw (3llr.center)--(2ll.center);
        \draw (3lml.center)--(2lm.center);
        \draw (3lmm.center)--(2lm.center);
        \draw (3lmr.center)--(2lm.center);
        \draw (3mll.center)--(2ml.center);
        \draw (3mlr.center)--(2ml.center);
        \draw (3mrl.center)--(2mr.center);
        \draw (3mrr.center)--(2mr.center);
        \draw (3lrl.center)--(2lr.center);
        \draw (3lrm.center)--(2lr.center);
        \draw (3lrr.center)--(2lr.center);
        \draw (3rll.center)--(2rl.center);
        \draw (3rlr.center)--(2rl.center);
        \draw (3rrl.center)--(2rr.center);
        \draw (3rrr.center)--(2rr.center);
    \end{tikzpicture}\begin{tikzpicture} [scale=.7]
        \node (0) at (0,0) {};
        \node (1l) at (-1,-1) {};
        \node (1r) at (1,-1) {};
        \node (2ll) at (-1.5,-2) {};
        \node (2lm) at (-1,-2) {};
        \node (2lr) at (-.5,-2) {};
        \node (2rl) at (.5,-2) {};
        \node (2rm) at (1,-2) {};
        \node (2rr) at (1.5,-2) {};
        \node (3lll) at (-1.6,-3) {};
        \node (3llm) at (-1.5,-3) {};
        \node (3llr) at (-1.4,-3) {};
        \node (3lml) at (-1.1,-3) {};
        \node (3lmr) at (-.9,-3) {};
        \node (3lrl) at (-.6,-3) {};
        \node (3lrm) at (-.5,-3) {};
        \node (3lrr) at (-.4,-3) {};
        \node (3rrr) at (1.6,-3) {};
        \node (3rrm) at (1.5,-3) {};
        \node (3rrl) at (1.4,-3) {};
        \node (3rmr) at (1.1,-3) {};
        \node (3rml) at (.9,-3) {};
        \node (3rlr) at (.6,-3) {};
        \node (3rlm) at (.5,-3) {};
        \node (3rll) at (.4,-3) {};
        \draw (2ll.center)--(1l.center)--(0.center)--(1r.center)--(2rr.center);
        \draw (2lm.center)--(1l.center);
        \draw (2lr.center)--(1l.center);
        \draw (2rl.center)--(1r.center);
        \draw (2rm.center)--(1r.center);
        \draw (3lll.center)--(2ll.center);
        \draw (3llm.center)--(2ll.center);
        \draw (3llr.center)--(2ll.center);
        \draw (3lml.center)--(2lm.center);
        \draw (3lmr.center)--(2lm.center);
        \draw (3lrl.center)--(2lr.center);
        \draw (3lrm.center)--(2lr.center);
        \draw (3lrr.center)--(2lr.center);
        \draw (3rrr.center)--(2rr.center);
        \draw (3rrm.center)--(2rr.center);
        \draw (3rrl.center)--(2rr.center);
        \draw (3rmr.center)--(2rm.center);
        \draw (3rml.center)--(2rm.center);
        \draw (3rlr.center)--(2rl.center);
        \draw (3rlm.center)--(2rl.center);
        \draw (3rll.center)--(2rl.center);
    \end{tikzpicture}}$$
        \caption{Construction of the humanoid to left-multiply with in order to swap two consecutive carets of the left tree.}
        \label{permutacarets2}
    \end{figure}
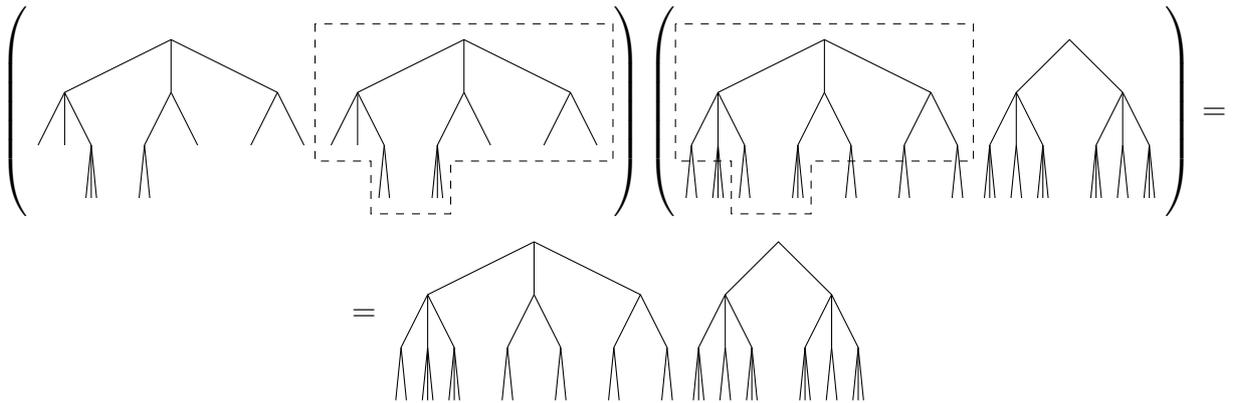

    \begin{lmm}\label{permutarcarets}
    Given a balanced tree-pair diagram, we can permute two consecutive carets on the bottom row of the left tree by left-multiplying by a humanoid. Similarly, we can permute two consecutive carets on the bottom row of the right tree by right-multiplying by a humanoid.
    \end{lmm}
    \begin{prf}
    In order to permute two consecutive carets on the bottom row of the left tree, we need to take as $T$ the left tree without the last row of carets, and match the added carets with the carets we want to swap. An example can be seen in \Cref{permutacarets2}. A similar construction can be done to permute carets on the right tree.
    \end{prf}

    Now, we are able to prove \Cref{humgen}:

    \begin{prf}[ of \Cref{humgen}]
    We need to write an arbitrary element as a product of humanoids. Equivalently, we will start with any element and we will left- and right-multiply by humanoids until we reach the identity. As every element admits a balanced tree representation, we will assume the tree-pair diagram of the given element is balanced. We will induct on the depth $d$ of its leaves. Since all tree-pair diagrams with all leaves at depth 0 or depth 1 are the identity, the base case is covered.

    It will be important for us to track the number of nodes in each tree that are at a depth one unit lower than the leaves. We call these nodes \textit{branches}. Notice that for each tree, if the number of leaves is $\ell$ and the number of branches is $b$, then $2b\leq\ell\leq3b$. The values of $b=n_2+n_3$ and $\ell=2n_2+3n_3$ are precisely determined by the number $n_2$ of $2$-carets and the number $n_3$ of $3$-carets. The converse is also true: the number of $2$-carets in the last row of the tree is exactly $n_2=3b-\ell$, and the number of $2$-carets is $n_3=\ell-2b$.

    If the number of branches of the left tree is the same as the number of branches of the right tree, then there are the same number of $2$-carets and $3$-carets in the last row, so by \Cref{permutarcarets} we can permute the ones on the left to match the ones on the right, making all carets in the last row superfluous. After those carets have been removed, we end up with an element of depth $d-1$, which by induction hypothesis can be written as a product of humanoids.

    If, instead, the number of branches does not coincide, let $b$ and $b'$ be the number of branches of the left tree and the right tree, respectively, and we have $b'\neq b$. We will only discuss the case $b'>b$ (the other case is analogous).

    We have the inequalities $2b\leq\ell\leq3b$ and $2b'\leq\ell\leq3b'$ which, combined with $b'\geq b+1$, give $2b+2\leq2b'\leq\ell\leq3b\leq3b'-3$. Therefore, there are at least two $3$-carets in the last row of the left tree, because $\ell-2b\geq(2b+2)-2b=2$. Now, we look for a $2$-caret of the left tree, excluding the ones in the last row. If there is none, all the carets hanging from nodes of depth $0,1,\dots,d-2$ are $3$-carets, so there will be $b=3^{d-1}$ branches on the left tree, and it would be impossible for the other tree to have more branches than that. Therefore, there will be some $2$-carets, and we can consider the bottommost one (or one of the bottommost ones).

    Unless this $2$-caret is on the second-to-last row of carets, its two children will be necessarily $3$-carets, so we substitute this collective by a $3$-caret with $2$-carets as children (see \Cref{2332}), obtaining a $2$-caret one row below. We repeat the process until there is a $2$-caret in the second-to-last row. Notice that the number of branches has remained constant during those transformations.

    \begin{figure}[H]
        $$\raisebox{-\height/2}{\begin{tikzpicture}
        \node (t) at (0,0) {};
        \node (l) at (-1,-1) {};
        \node (r) at (1,-1) {};
        \node (ll) at (-1.25,-2) {};
        \node (lm) at (-1,-2) {};
        \node (lr) at (-0.75,-2) {};
        \node (rl) at (0.75,-2) {};
        \node (rm) at (1,-2) {};
        \node (rr) at (1.25,-2) {};
        \draw (ll.center)--(l.center)--(t.center)--(r.center)--(rr.center);
        \draw (l.center)--(lr.center);
        \draw (lm.center)--(l.center);
        \draw (r.center)--(rm.center);
        \draw (rl.center)--(r.center);
        \node (T) at (3,0) {};
        \node (L) at (2,-1) {};
        \node (M) at (3,-1) {};
        \node (R) at (4,-1) {};
        \node (LL) at (1.75,-2) {};
        \node (LR) at (2.25,-2) {};
        \node (RL) at (3.75,-2) {};
        \node (RR) at (4.25,-2) {};
        \node (ML) at (2.75,-2) {};
        \node (MR) at (3.25,-2) {};
        \draw (LL.center)--(L.center)--(T.center)--(R.center)--(RR.center);
        \draw (L.center)--(LR.center);
        \draw (ML.center)--(M.center);
        \draw (M.center)--(T.center);
        \draw (M.center)--(MR.center);
        \draw (RL.center)--(R.center);
        \draw[<->] (1,-.5)--(2,-.5);
    \end{tikzpicture}}$$
    \caption{Transformation of a tree that does not alter the element.}
    \label{2332}
    \end{figure}
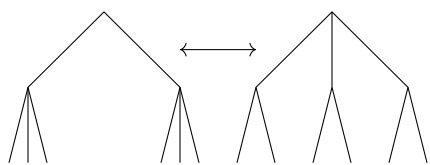

    Now that we have a $2$-caret in the second-to-last row, we can permute the last row until two $3$-carets are located below the $2$-caret. Performing once more the aforementioned transformation changes the number of branches of the first tree, increasing it by one and therefore closing in the gap between $b$ and $b'$ by one unit. By repeatedly doing this, we eventually make $b'=b$ and we are back to the previous solved case.
    \end{prf}

    This proves that humanoids generate $F\left(\frac32\right)$. Before we go further, we need to address why we called them this way. Notice that as the same tree $T$ is used on both sides, there is potentially a lot of cancellation. The only carets that do not disappear are the ones that embrace at least one of the bottommost carets.

    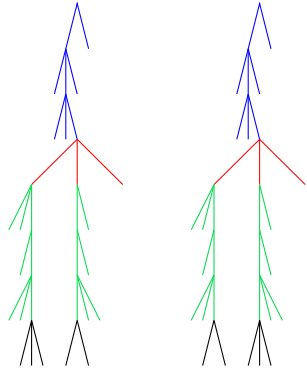
\begin{figure}[H]
    \centering
    \raisebox{-\height/2}{\begin{tikzpicture}[scale=.6]
        \node (0) at (0,0) {};
        \node (1l) at (-.25,-1) {};
        \node (1r) at (.25,-1) {};
        \node (2l) at (-.5,-2) {};
        \node (2m) at (-.25,-2) {};
        \node (2r) at (0,-2) {};
        \node (3l) at (-.5,-3) {};
        \node (3m) at (-.25,-3) {};
        \node (3r) at (0,-3) {};
        \node (4l) at (-1,-4) {};
        \node (4m) at (0,-4) {};
        \node (4r) at (1,-4) {};
        \node (5L) at (-1.5,-5) {};
        \node (5M) at (-1.25,-5) {};
        \node (5R) at (-1,-5) {};
        \node (6L) at (-1.25,-6) {};
        \node (6R) at (-1,-6) {};
        \node (7L) at (-1.5,-7) {};
        \node (7M) at (-1.25,-7) {};
        \node (7R) at (-1,-7) {};
        \node (8L) at (-1.25,-8) {};
        \node (8M) at (-1,-8) {};
        \node (8R) at (-0.75,-8) {};
        \node (5l) at (0,-5) {};
        \node (5r) at (.25,-5) {};
        \node (6l) at (0,-6) {};
        \node (6r) at (.25,-6) {};
        \node (7l) at (0,-7) {};
        \node (7m) at (.25,-7) {};
        \node (7r) at (.5,-7) {};
        \node (8l) at (-.25,-8) {};
        \node (8r) at (.25,-8) {};
        \draw[color=red] (4l.center)--(3r.center);
        \draw[color=red] (4m.center)--(3r.center);
        \draw[color=red] (4r.center)--(3r.center);
        \draw[color=blue] (3r.center)--(2m.center);
        \draw[color=blue] (3m.center)--(2m.center)--(1l.center);
        \draw[color=blue] (3l.center)--(2m.center);
        \draw[color=blue] (2l.center)--(1l.center)--(0.center)--(1r.center);
        \draw[color=blue] (2r.center)--(1l.center);
        \draw[color=malaquita] (7R.center)--(6R.center)--(5R.center)--(4l.center);
        \draw[color=malaquita] (7M.center)--(6R.center);
        \draw[color=malaquita] (7L.center)--(6R.center);
        \draw[color=malaquita] (6L.center)--(5R.center);
        \draw[color=malaquita] (5M.center)--(4l.center);
        \draw[color=malaquita] (5L.center)--(4l.center);
        \draw[color=malaquita] (5r.center)--(4m.center);
        \draw[color=malaquita] (6r.center)--(5l.center);
        \draw[color=malaquita] (7m.center)--(6l.center);
        \draw[color=malaquita] (7r.center)--(6l.center);
        \draw[color=malaquita] (7l.center)--(6l.center)--(5l.center)--(4m.center);
        \draw (8l.center)--(7l.center);
        \draw (8r.center)--(7l.center);
        \draw (8L.center)--(7R.center);
        \draw (8M.center)--(7R.center);
        \draw (8R.center)--(7R.center);
        \node (0') at (4,0) {};
        \node (1l') at (4-.25,-1) {};
        \node (1r') at (4.25,-1) {};
        \node (2l') at (4-.5,-2) {};
        \node (2m') at (4-.25,-2) {};
        \node (2r') at (4,-2) {};
        \node (3l') at (4-.5,-3) {};
        \node (3m') at (4-.25,-3) {};
        \node (3r') at (4,-3) {};
        \node (4l') at (4-1,-4) {};
        \node (4m') at (4,-4) {};
        \node (4r') at (4+1,-4) {};
        \node (5L') at (4-1.5,-5) {};
        \node (5M') at (4-1.25,-5) {};
        \node (5R') at (4-1,-5) {};
        \node (6L') at (4-1.25,-6) {};
        \node (6R') at (4-1,-6) {};
        \node (7L') at (4-1.5,-7) {};
        \node (7M') at (4-1.25,-7) {};
        \node (7R') at (4-1,-7) {};
        \node (8L') at (4-1.25,-8) {};
        \node (8R') at (4-0.75,-8) {};
        \node (5l') at (4,-5) {};
        \node (5r') at (4.25,-5) {};
        \node (6l') at (4,-6) {};
        \node (6r') at (4.25,-6) {};
        \node (7l') at (4,-7) {};
        \node (7m') at (4.25,-7) {};
        \node (7r') at (4.5,-7) {};
        \node (8l') at (4-.25,-8) {};
        \node (8m') at (4,-8) {};
        \node (8r') at (4.25,-8) {};
        \draw[color=red] (4l'.center)--(3r'.center);
        \draw[color=red] (4m'.center)--(3r'.center);
        \draw[color=red] (4r'.center)--(3r'.center);
        \draw[color=blue] (3r'.center)--(2m'.center);
        \draw[color=blue] (3m'.center)--(2m'.center)--(1l'.center);
        \draw[color=blue] (3l'.center)--(2m'.center);
        \draw[color=blue] (2l'.center)--(1l'.center)--(0'.center)--(1r'.center);
        \draw[color=blue] (2r'.center)--(1l'.center);
        \draw[color=malaquita] (7R'.center)--(6R'.center)--(5R'.center)--(4l'.center);
        \draw[color=malaquita] (7M'.center)--(6R'.center);
        \draw[color=malaquita] (7L'.center)--(6R'.center);
        \draw[color=malaquita] (6L'.center)--(5R'.center);
        \draw[color=malaquita] (5M'.center)--(4l'.center);
        \draw[color=malaquita] (5L'.center)--(4l'.center);
        \draw[color=malaquita] (5r'.center)--(4m'.center);
        \draw[color=malaquita] (6r'.center)--(5l'.center);
        \draw[color=malaquita] (7m'.center)--(6l'.center);
        \draw[color=malaquita] (7r'.center)--(6l'.center);
        \draw[color=malaquita] (7l'.center)--(6l'.center)--(5l'.center)--(4m'.center);
        \draw (8l'.center)--(7l'.center);
        \draw (8m'.center)--(7l'.center);
        \draw (8r'.center)--(7l'.center);
        \draw (8L'.center)--(7R'.center);
        \draw (8R'.center)--(7R'.center);
    \end{tikzpicture}}
    \caption{Example of a humanoid element after the cancellation.}
    \end{figure}

    \begin{dfn}
    We define the \emph{hip} of a humanoid tree-pair diagram (red in the example above) as the bottommost caret of $T$ that embraces the two consecutive leaves where the added carets hang from. The carets above the hip are called the \emph{torso} (blue in the example above), while the carets of $T$ embraced by the hip will be called the \emph{legs} (green in the example above). Finally, the two added carets will be the \emph{feet} (black in the example above).
    \end{dfn}

    It is because of the resemblance with the human body that we call these elements humanoid. Note that every humanoid has a hip and two feet, but it may not have a torso, and it may not have legs. If we need to refer to a particular leg or foot, we will do so as if we are seeing the humanoid from behind, so that its left leg/foot is on our left and the right one is on our right.

    \section{Serpents}
    We continue with the next family, which gets the name from the resemblance to snake-like creatures:\\

    \begin{dfn}
    A humanoid is called a \textit{serpent} if it has no legs. That is, a humanoid whose feet hang directly from the hip.\\
    \end{dfn}

    \begin{prp}\label{serpgen}
    Serpents generate $F\left(\frac32\right)$.
    \end{prp}

    \begin{prf} It is sufficient to generate every humanoid. For the sake of contradiction, suppose that there are humanoids that cannot be written in terms of serpents. Then, take all humanoid tree-pair diagrams of all of those elements and choose the ones that have the minimum depth, and of those, one that has the shortest legs. Since the element is not a serpent, the length of the legs is at least 1.

    We claim that the legs need to be symmetrical. By that, we mean that the type of carets we visit when going from the hip to the feet is the same for both legs at every step. If it is not the case, at one moment we will visit a $3$-caret on one leg and a $2$-caret on the other one. In that case, we construct the element (see \Cref{amputacames}), whose right tree is obtained by taking the left tree of the tree-pair diagram but with its legs cut off, so that the asymmetric part of the leg becomes the feet. The left tree of this new element is the same but with its feet interchanged. This new element has strictly less depth, so can be constructed using serpents by hypothesis.\\

    \begin{figure}[H]
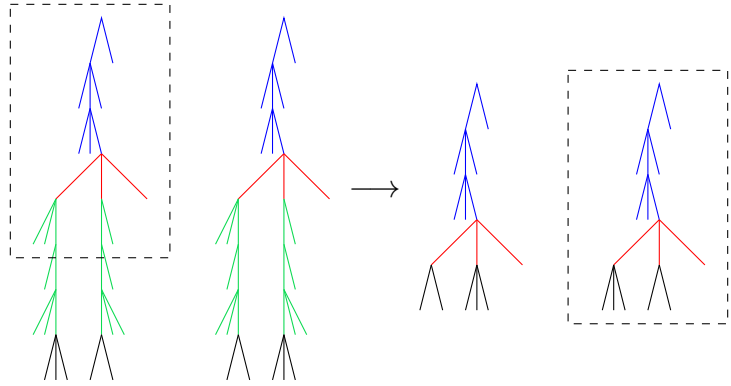

        $$\raisebox{-\height/2}{
}$$
    \caption{How to construct the element we will conjugate by in order to lower the hip.}
        \label{amputacames}
    \end{figure}

    We conjugate by this new element. The result is that the hip lowers to the current part of the leg (see \Cref{baixahip}), so the legs shorten and the result is also in the span of the serpents. Therefore, the original element itself also belongs to the span.

    \begin{figure}[H]
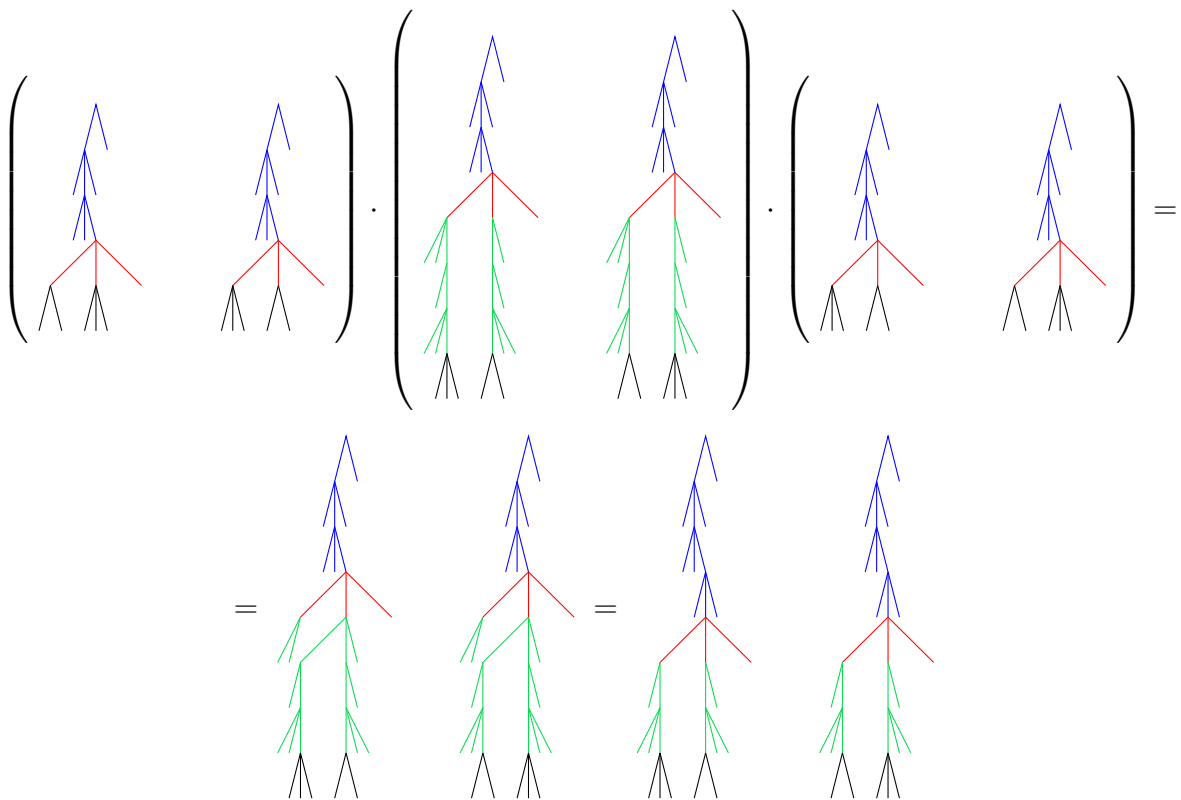

    $$\left(\raisebox{-\height/2}{%
}$$
    \caption{Conjugation lowers the hip.}
    \label{baixahip}
    \end{figure}

    So, we conclude that all the alleged humanoids not generated by serpents that have minimal depth and lowest hip, must have symmetrical legs. Moreover, we now say that the legs can have at most one caret each. Otherwise, similarly to before, we conjugate by another element, which we construct as follows.\\

    Take the left tree and remove the right leg and foot. The result is a sequence of carets, each one hanging from the previous one. Remove the last one (which used to be the left foot) and consider the caret it was hanging from (which used to be the last part of the leg). Call that caret $A$. Since the leg was at least two units long, caret $A$ hangs from the right stem of another caret (which was also part of the leg), called caret $B$. Now, add a caret $C$ to the left of $A$, of type opposite to $A$, which hangs from a contiguous leaf of $B$. The resulting tree is the right tree of a serpent, which has $B$ as the hip and $C$ and $A$ as its feet. The left tree of this serpent is obtained by interchanging $A$ and $C$. Refer to \Cref{serp} for an example.

    \begin{figure}[H]
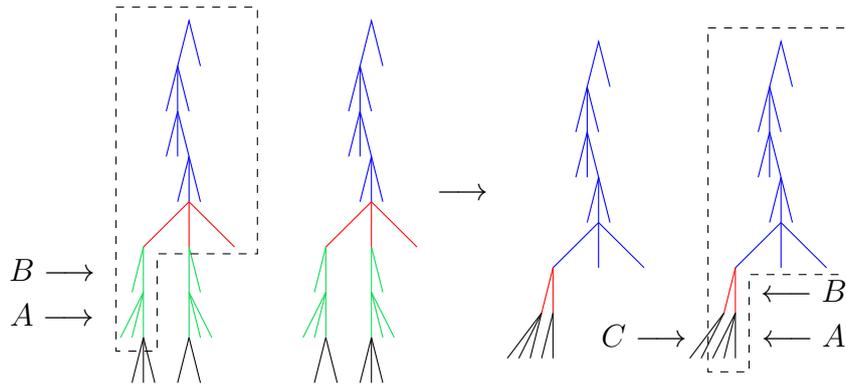

        $$\raisebox{-\height/2}{
}$$
        \caption{How to construct the element we will conjugate by in order to make legs asymmetrical.}
        \label{serp}
    \end{figure}

    We constructed a serpent and we conjugate our element with it. The tree-pair diagram we end up with (see \Cref{escurçacames}) is the same as the one we started with, except that the last left-leg caret of both trees changed to the other type: it has the same depth and leg length as before, but asymmetric legs, so it is generated by humanoids.

    \begin{figure}[H]
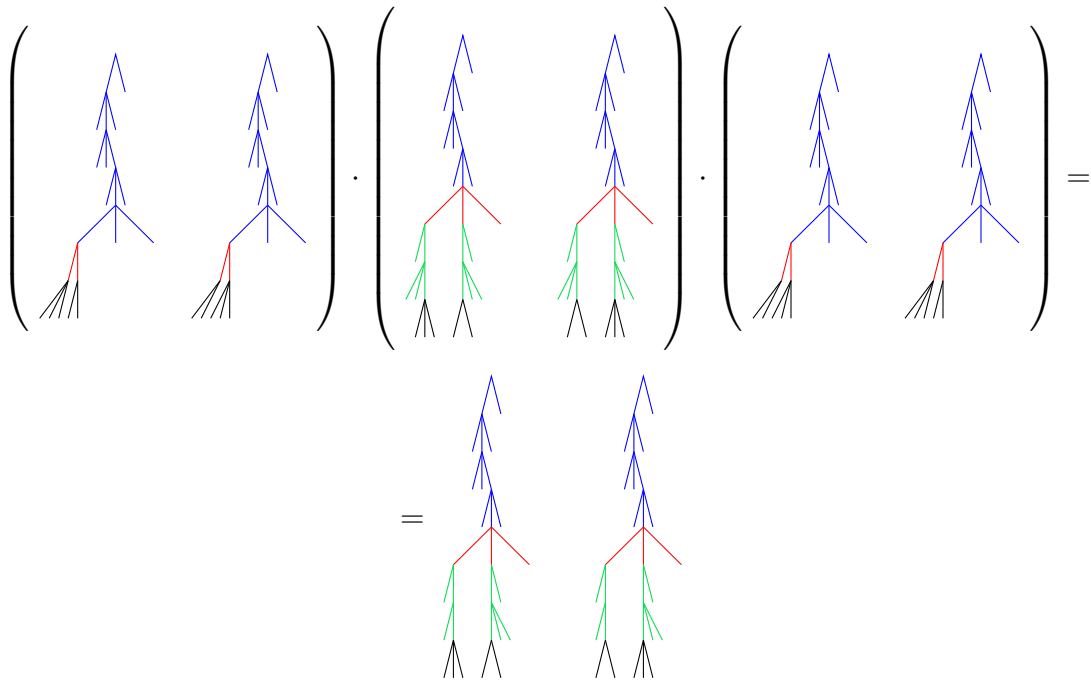

        $$\left(\raisebox{-\height/2}{%
}$$
        \caption{Conjugation makes legs asymmetrical.}
        \label{escurçacames}
    \end{figure}

    What we have seen so far is that if there is any humanoid that may not be obtained by combining serpents, then there is also one with symmetric legs of length exactly $1$. However all of those are generated by serpents, thanks to the following equalities, which the reader can verify (the respective equalities for the inverses can easily be deduced from the given ones). In each equality, the three dots, wherever they appear, represent the same arrangement of carets, which constitute the torso of the tree-pair diagram in the left-hand side member of the equality.
    \end{prf}

    \begin{figure}[H]
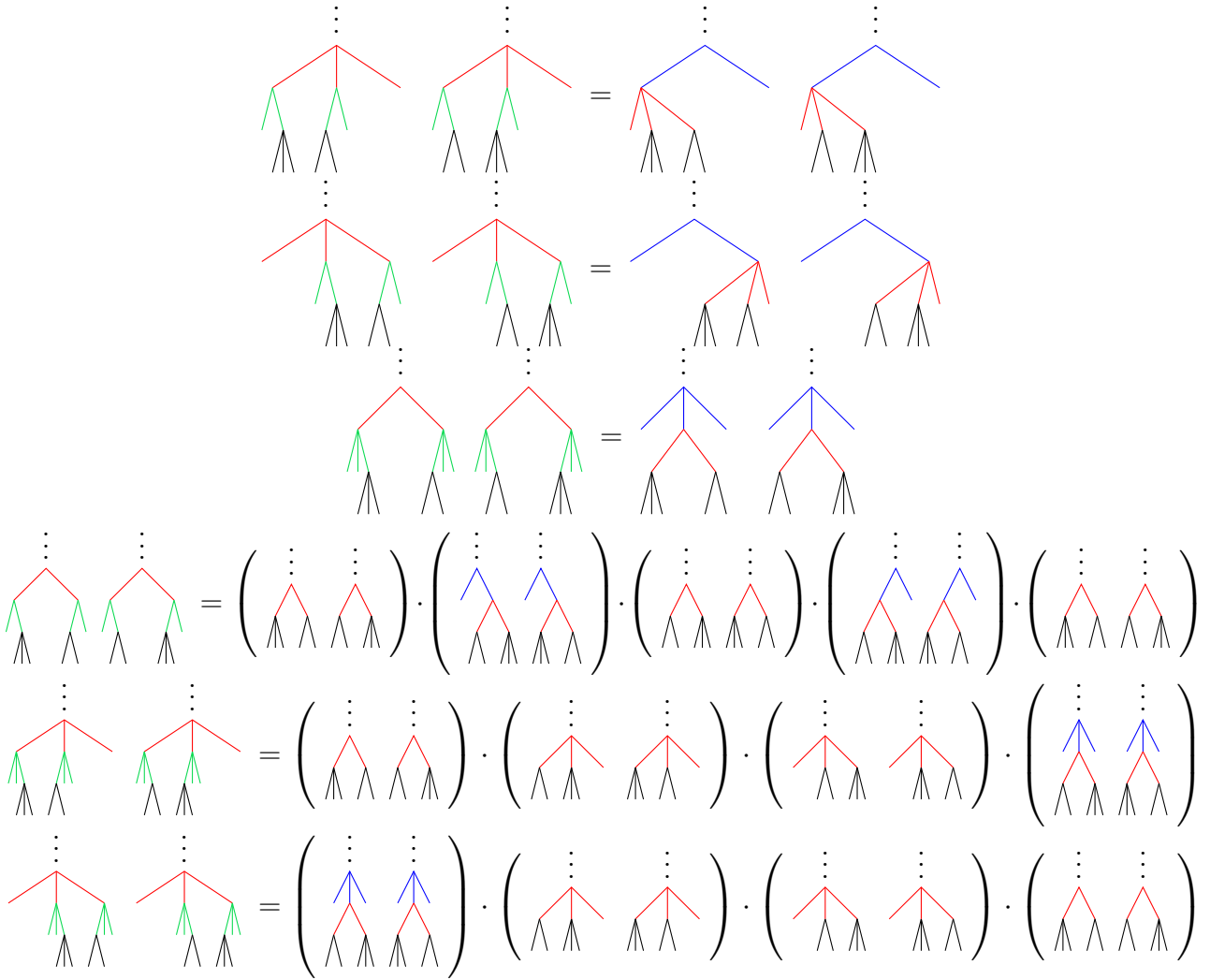

     $$\raisebox{-\height/2}{
}\right)$$
    \caption{Decomposition of one-long-symmetrically-legged humanoids into serpents.}
    \end{figure}

    Note that the first three equalities are not products: we can pass from the left-hand side to the right-hand side by the transformation shown in \Cref{2332}. On the first two, a $2$-caret needs to be added on both trees before the transformation can be applied. For the other three relations, we need to repeatedly find common subdivisions in order to perform the multiplications. It is tedious but not difficult, and this is why it was left for the reader to check.
    \newpage
    \section{Centipedes}

    We continue our parade of generating families:

    \begin{dfn} A \textit{centipede} is a serpent whose torso and hip carets are only allowed to be $3$-carets.
    \end{dfn}

    The name comes from the fact that a centipede is a snake-like creature with a great number of limbs.

    \begin{prp}\label{centgenera} Centipedes generate $F\left(\frac32\right)$.
    \end{prp}

\begin{figure}[H]
    $$\raisebox{-\height/2}{\begin{tikzpicture}[scale=.6]
        \node (dots) at (0,.3) {$\smash{\vdots}$};
        \node (0) at (0,0) {};
        \node (1l) at (-.5,-1) {};
        \node (1r) at (.5,-1) {};
        \node (2l) at (-.75,-2) {};
        \node (2m) at (-.5,-2) {};
        \node (2r) at (-.25,-2) {};
        \node (2L) at (.25,-2) {};
        \node (2R) at (.75,-2) {};
        \draw[color=red] (1l.center)--(0.center)--(1r.center);
        \draw (2l.center)--(1l.center);
        \draw (1l.center)--(2r.center);
        \draw (2m.center)--(1l.center);
        \draw (2L.center)--(1r.center);
        \draw (1r.center)--(2R.center);
        \node (Dots) at (2,.3) {$\smash{\vdots}$};
        \node (0') at (2,0) {};
        \node (1l') at (1.5,-1) {};
        \node (1r') at (2.5,-1) {};
        \node (2l') at (1.25,-2) {};
        \node (2r') at (1.75,-2) {};
        \node (2L') at (2.25,-2) {};
        \node (2M') at (2.5,-2) {};
        \node (2R') at (2.75,-2) {};
        \draw[color=red] (1l'.center)--(0'.center)--(1r'.center);
        \draw (2l'.center)--(1l'.center);
        \draw (1l'.center)--(2r'.center);
        \draw (2L'.center)--(1r'.center);
        \draw (2M'.center)--(1r'.center);
        \draw (1r'.center)--(2R'.center);
    \end{tikzpicture}}=\left(\raisebox{-\height/2}{\begin{tikzpicture}[scale=.6]
        \node (dots) at (.5,.3) {$\smash{\vdots}$};
        \node (0) at (.5,0) {};
        \node (1l) at (-.5,-1) {};
        \node (1r) at (.5,-1) {};
        \node (1rr) at (1.5,-1) {};
        \node (2l) at (-.75,-2) {};
        \node (2m) at (-.5,-2) {};
        \node (2r) at (-.25,-2) {};
        \node (2L) at (.25,-2) {};
        \node (2R) at (.75,-2) {};
        \draw[color=red] (1l.center)--(0.center)--(1rr.center);
        \draw (2l.center)--(1l.center);
        \draw (1l.center)--(2r.center);
        \draw (2m.center)--(1l.center);
        \draw (2L.center)--(1r.center);
        \draw (1r.center)--(2R.center);
        \draw[color=red] (1r.center)--(0.center);
        \node (Dots) at (3.5,.3) {$\smash{\vdots}$};
        \node (0') at (1+2.5,0) {};
        \node (1l') at (1+1.5,-1) {};
        \node (1r') at (1+2.5,-1) {};
        \node (1rr') at (2+2.5,-1) {};
        \node (2l') at (1+1.25,-2) {};
        \node (2r') at (1+1.75,-2) {};
        \node (2L') at (1+2.25,-2) {};
        \node (2M') at (3.5,-2) {};
        \node (2R') at (1+2.75,-2) {};
        \draw[color=red] (1l'.center)--(0'.center)--(1rr'.center);
        \draw (2l'.center)--(1l'.center);
        \draw (1l'.center)--(2r'.center);
        \draw (2L'.center)--(1r'.center);
        \draw (2M'.center)--(1r'.center);
        \draw (1r'.center)--(2R'.center);
        \draw[color=red] (1r'.center)--(0'.center);
    \end{tikzpicture}}\right)\cdot\left(\raisebox{-\height/2}{\begin{tikzpicture}[scale=.6]
        \node (dots) at (.5,.3) {$\smash{\vdots}$};
        \node (0) at (.5,0) {};
        \node (1l) at (-.5,-1) {};
        \node (1r) at (.5,-1) {};
        \node (1rr) at (1.5,-1) {};
        \node (2l) at (1-.75,-2) {};
        \node (2m) at (2-.5-1,-2) {};
        \node (2r) at (1-.25,-2) {};
        \node (2L) at (1.25,-2) {};
        \node (2R) at (1.75,-2) {};
        \draw[color=red] (1l.center)--(0.center)--(1rr.center);
        \draw (2l.center)--(1r.center);
        \draw (1r.center)--(2r.center);
        \draw (2m.center)--(1r.center);
        \draw (2L.center)--(1rr.center);
        \draw (1rr.center)--(2R.center);
        \draw[color=red] (1r.center)--(0.center);
        \node (Dots) at (3.5,.3) {$\smash{\vdots}$};
        \node (0') at (1+2.5,0) {};
        \node (1l') at (1+1.5,-1) {};
        \node (1r') at (1+2.5,-1) {};
        \node (1rr') at (2+2.5,-1) {};
        \node (2l') at (2+1.25,-2) {};
        \node (2r') at (2+1.75,-2) {};
        \node (2L') at (2+2.25,-2) {};
        \node (2M') at (3.5+1,-2) {};
        \node (2R') at (2+2.75,-2) {};
        \draw[color=red] (1l'.center)--(0'.center)--(1rr'.center);
        \draw (2l'.center)--(1r'.center);
        \draw (1r'.center)--(2r'.center);
        \draw (2L'.center)--(1rr'.center);
        \draw (2M'.center)--(1rr'.center);
        \draw (1rr'.center)--(2R'.center);
        \draw[color=red] (1r'.center)--(0'.center);
    \end{tikzpicture}}\right)$$
    \caption{Decomposing a serpent with a $2$-caret hip by serpents with a $3$-caret hip.}
    \label{hip2}
    \end{figure}
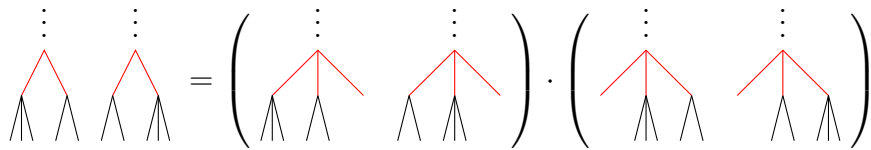

    In order to prove \Cref{centgenera}, we need to generate every serpent, and we first focus on the hip:

    \begin{lmm}\label{hip} Every serpent can be written as a product of at most two serpents whose hip is a $3$-caret.
    \end{lmm}
    \begin{prf}
        If the serpent has a $3$-caret hip, there is nothing left to prove. If the hip is a $2$-caret, we decompose the element as in \Cref{hip2}.
    \end{prf}

    \begin{prf}[ of \Cref{centgenera}] Because of \Cref{hip}, we will only work with serpents with a $3$-caret hip. We will perform an induction on the number of $2$-carets in the torso. If no $2$-caret is present in the torso, the element is a centipede and we are done. So, we suppose there are some. If the bottommost $2$-caret in the torso is not right above the hip, we can add a superfluous $3$-caret on its free stem (on both trees) and perform the transformation in \Cref{2332}, after which we will be able to remove two superfluous pairs of $2$-carets and end up with a serpent with the same number of $2$-carets in its torso as before, but the bottommost $2$-caret will be a unit lower.

    We can repeat the procedure until the $2$-caret is just above the hip. Then, we can perform one of the following decompositions, depending on the relative position of the $2$-caret and the hip. To save space, we used elements with a $2$-caret hip, which ought to be further decomposed according to \Cref{hip}.

    \begin{figure}[H]
    $$\raisebox{-\height/2}{\begin{tikzpicture}[scale=.6]
        \node (dots) at (2,.3) {$\smash{\vdots}$};
        \node (0) at (2,0) {};
        \node (1l) at (3,-1) {};
        \node (1r) at (1,-1) {};
        \node (2l) at (2,-2) {};
        \node (2r) at (0,-2) {};
        \node (2L) at (1,-2) {};
        \node (3L) at (.75,-3) {};
        \node (3R) at (1.25,-3) {};
        \node (3r) at (.25,-3) {};
        \node (3m) at (0,-3) {};
        \node (3l) at (-.25,-3) {};
        \draw[color=blue] (1l.center)--(0.center)--(1r.center);
        \draw[color=red] (2l.center)--(1r.center);
        \draw[color=red] (1r.center)--(2r.center);
        \draw[color=red] (2L.center)--(1r.center);
        \draw (3l.center)--(2r.center);
        \draw (3m.center)--(2r.center);
        \draw (2r.center)--(3r.center);
        \draw (3L.center)--(2L.center);
        \draw (2L.center)--(3R.center);
        \node (Dots) at (6,.3) {$\smash{\vdots}$};
        \node (0') at (6,0) {};
        \node (1l') at (7,-1) {};
        \node (1r') at (5,-1) {};
        \node (2l') at (6,-2) {};
        \node (2r') at (4,-2) {};
        \node (2L') at (5,-2) {};
        \node (3L') at (4.75,-3) {};
        \node (3M') at (5,-3) {};
        \node (3R') at (5.25,-3) {};
        \node (3l') at (4.25,-3) {};
        \node (3r') at (3.75,-3) {};
        \draw[color=blue] (1l'.center)--(0'.center)--(1r'.center);
        \draw[color=red] (2l'.center)--(1r'.center);
        \draw[color=red] (1r'.center)--(2r'.center);
        \draw[color=red] (2L'.center)--(1r'.center);
        \draw (3l'.center)--(2r'.center);
        \draw (3M'.center)--(2L'.center);
        \draw (2r'.center)--(3r'.center);
        \draw (3L'.center)--(2L'.center);
        \draw (2L'.center)--(3R'.center);
    \end{tikzpicture}}=\raisebox{-\height/2}{\begin{tikzpicture}[scale=.6]
        \node (dots) at (3,.3) {$\smash{\vdots}$};
        \node (0) at (3,0) {};
        \node (1l) at (4.5,-1) {};
        \node (1m) at (3,-1) {};
        \node (1r) at (1.5,-1) {};
        \node (2l) at (1,-2) {};
        \node (2r) at (2,-2) {};
        \node (3L) at (1.75,-3) {};
        \node (3R) at (2.25,-3) {};
        \node (3r) at (.75,-3) {};
        \node (3m) at (1,-3) {};
        \node (3l) at (1.25,-3) {};
        \draw[color=blue] (1m.center)--(0.center);
        \draw[color=blue] (1l.center)--(0.center)--(1r.center);
        \draw[color=red] (2l.center)--(1r.center);
        \draw[color=red] (1r.center)--(2r.center);
        \draw[color=red] (2L.center)--(1r.center);
        \draw (3l.center)--(2l.center);
        \draw (3m.center)--(2l.center);
        \draw (2l.center)--(3r.center);
        \draw (3L.center)--(2r.center);
        \draw (2r.center)--(3R.center);
        \node (dots) at (7,.3) {$\smash{\vdots}$};
        \node (0') at (7,0) {};
        \node (1l') at (8.5,-1) {};
        \node (1m') at (7,-1) {};
        \node (1r') at (5.5,-1) {};
        \node (2r') at (5,-2) {};
        \node (2L') at (6,-2) {};
        \node (3L') at (5.75,-3) {};
        \node (3M') at (6,-3) {};
        \node (3R') at (6.25,-3) {};
        \node (3l') at (5.25,-3) {};
        \node (3r') at (4.75,-3) {};
        \draw[color=blue] (1m'.center)--(0'.center);
        \draw[color=blue] (1l'.center)--(0'.center)--(1r'.center);
        \draw[color=red] (1r'.center)--(2r'.center);
        \draw[color=red] (2L'.center)--(1r'.center);
        \draw (3l'.center)--(2r'.center);
        \draw (3M'.center)--(2L'.center);
        \draw (2r'.center)--(3r'.center);
        \draw (3L'.center)--(2L'.center);
        \draw (2L'.center)--(3R'.center);
    \end{tikzpicture}}$$
    $$\raisebox{-\height/2}{\begin{tikzpicture}[scale=.6]
        \node (dots) at (0,.3) {$\smash{\vdots}$};
        \node (0) at (0,0) {};
        \node (1l) at (-1,-1) {};
        \node (1r) at (1,-1) {};
        \node (2l) at (0,-2) {};
        \node (2r) at (1,-2) {};
        \node (2L) at (2,-2) {};
        \node (3L) at (1.75,-3) {};
        \node (3R) at (2.25,-3) {};
        \node (3r) at (1.25,-3) {};
        \node (3m) at (1,-3) {};
        \node (3l) at (.75,-3) {};
        \draw[color=blue] (1l.center)--(0.center)--(1r.center);
        \draw[color=red] (2l.center)--(1r.center);
        \draw[color=red] (1r.center)--(2r.center);
        \draw[color=red] (2L.center)--(1r.center);
        \draw (3l.center)--(2r.center);
        \draw (3m.center)--(2r.center);
        \draw (2r.center)--(3r.center);
        \draw (3L.center)--(2L.center);
        \draw (2L.center)--(3R.center);
        \node (dots) at (4,.3) {$\smash{\vdots}$};
        \node (0') at (4,0) {};
        \node (1l') at (3,-1) {};
        \node (1r') at (5,-1) {};
        \node (2l') at (4,-2) {};
        \node (2r') at (5,-2) {};
        \node (2L') at (6,-2) {};
        \node (3L') at (5.75,-3) {};
        \node (3M') at (6,-3) {};
        \node (3R') at (6.25,-3) {};
        \node (3l') at (5.25,-3) {};
        \node (3r') at (4.75,-3) {};
        \draw[color=blue] (1l'.center)--(0'.center)--(1r'.center);
        \draw[color=red] (2l'.center)--(1r'.center);
        \draw[color=red] (1r'.center)--(2r'.center);
        \draw[color=red] (2L'.center)--(1r'.center);
        \draw (3l'.center)--(2r'.center);
        \draw (3M'.center)--(2L'.center);
        \draw (2r'.center)--(3r'.center);
        \draw (3L'.center)--(2L'.center);
        \draw (2L'.center)--(3R'.center);
    \end{tikzpicture}}=\raisebox{-\height/2}{\begin{tikzpicture}[scale=.6]
        \node (dots) at (0,.3) {$\smash{\vdots}$};
        \node (0) at (0,0) {};
        \node (1l) at (-1.5,-1) {};
        \node (1m) at (0,-1) {};
        \node (1r) at (1.5,-1) {};
        \node (2l) at (1,-2) {};
        \node (2r) at (2,-2) {};
        \node (3L) at (1.75,-3) {};
        \node (3R) at (2.25,-3) {};
        \node (3r) at (.75,-3) {};
        \node (3m) at (1,-3) {};
        \node (3l) at (1.25,-3) {};
        \draw[color=blue] (1m.center)--(0.center);
        \draw[color=blue] (1l.center)--(0.center)--(1r.center);
        \draw[color=red] (2l.center)--(1r.center);
        \draw[color=red] (1r.center)--(2r.center);
        \draw[color=red] (2L.center)--(1r.center);
        \draw (3l.center)--(2l.center);
        \draw (3m.center)--(2l.center);
        \draw (2l.center)--(3r.center);
        \draw (3L.center)--(2r.center);
        \draw (2r.center)--(3R.center);
        \node (Dots) at (4,.3) {$\smash{\vdots}$};
        \node (0') at (4,0) {};
        \node (1l') at (2.5,-1) {};
        \node (1m') at (4,-1) {};
        \node (1r') at (5.5,-1) {};
        \node (2r') at (5,-2) {};
        \node (2L') at (6,-2) {};
        \node (3L') at (5.75,-3) {};
        \node (3M') at (6,-3) {};
        \node (3R') at (6.25,-3) {};
        \node (3l') at (5.25,-3) {};
        \node (3r') at (4.75,-3) {};
        \draw[color=blue] (1m'.center)--(0'.center);
        \draw[color=blue] (1l'.center)--(0'.center)--(1r'.center);
        \draw[color=red] (1r'.center)--(2r'.center);
        \draw[color=red] (2L'.center)--(1r'.center);
        \draw (3l'.center)--(2r'.center);
        \draw (3M'.center)--(2L'.center);
        \draw (2r'.center)--(3r'.center);
        \draw (3L'.center)--(2L'.center);
        \draw (2L'.center)--(3R'.center);
    \end{tikzpicture}}$$
    $$\raisebox{-\height/2}{\begin{tikzpicture}[scale=.6]
        \node (dots) at (2,.3) {$\smash{\vdots}$};
        \node (0) at (2,0) {};
        \node (1l) at (3,-1) {};
        \node (1r) at (1,-1) {};
        \node (2l) at (0,-2) {};
        \node (2r) at (1,-2) {};
        \node (2L) at (2,-2) {};
        \node (3L) at (1.75,-3) {};
        \node (3R) at (2.25,-3) {};
        \node (3r) at (1.25,-3) {};
        \node (3m) at (1,-3) {};
        \node (3l) at (.75,-3) {};
        \draw[color=blue] (1l.center)--(0.center)--(1r.center);
        \draw[color=red] (2l.center)--(1r.center);
        \draw[color=red] (1r.center)--(2r.center);
        \draw[color=red] (2L.center)--(1r.center);
        \draw (3l.center)--(2r.center);
        \draw (3m.center)--(2r.center);
        \draw (2r.center)--(3r.center);
        \draw (3L.center)--(2L.center);
        \draw (2L.center)--(3R.center);
        \node (dots) at (6,.3) {$\smash{\vdots}$};
        \node (0') at (6,0) {};
        \node (1l') at (7,-1) {};
        \node (1r') at (5,-1) {};
        \node (2l') at (4,-2) {};
        \node (2r') at (5,-2) {};
        \node (2L') at (6,-2) {};
        \node (3L') at (5.75,-3) {};
        \node (3M') at (6,-3) {};
        \node (3R') at (6.25,-3) {};
        \node (3l') at (5.25,-3) {};
        \node (3r') at (4.75,-3) {};
        \draw[color=blue] (1l'.center)--(0'.center)--(1r'.center);
        \draw[color=red] (2l'.center)--(1r'.center);
        \draw[color=red] (1r'.center)--(2r'.center);
        \draw[color=red] (2L'.center)--(1r'.center);
        \draw (3l'.center)--(2r'.center);
        \draw (3M'.center)--(2L'.center);
        \draw (2r'.center)--(3r'.center);
        \draw (3L'.center)--(2L'.center);
        \draw (2L'.center)--(3R'.center);
    \end{tikzpicture}}=\left(\raisebox{-\height/2}{\begin{tikzpicture}[scale=.6]
        \node (dots) at (0,.3) {$\smash{\vdots}$};
        \node (0) at (0,0) {};
        \node (1l) at (-.5,-1) {};
        \node (1r) at (.5,-1) {};
        \node (2l) at (-.75,-2) {};
        \node (2m) at (.5,-2) {};
        \node (2r) at (-.25,-2) {};
        \node (2L) at (.25,-2) {};
        \node (2R) at (.75,-2) {};
        \draw[color=red] (1l.center)--(0.center)--(1r.center);
        \draw (2l.center)--(1l.center);
        \draw (1l.center)--(2r.center);
        \draw (2m.center)--(1r.center);
        \draw (2L.center)--(1r.center);
        \draw (1r.center)--(2R.center);
        \node (dots) at (2,.3) {$\smash{\vdots}$};
        \node (0') at (2,0) {};
        \node (1l') at (1.5,-1) {};
        \node (1r') at (2.5,-1) {};
        \node (2l') at (1.25,-2) {};
        \node (2r') at (1.75,-2) {};
        \node (2L') at (2.25,-2) {};
        \node (2M') at (1.5,-2) {};
        \node (2R') at (2.75,-2) {};
        \draw[color=red] (1l'.center)--(0'.center)--(1r'.center);
        \draw (2l'.center)--(1l'.center);
        \draw (1l'.center)--(2r'.center);
        \draw (2L'.center)--(1r'.center);
        \draw (2M'.center)--(1l'.center);
        \draw (1r'.center)--(2R'.center);
    \end{tikzpicture}}\right)\cdot\left(\raisebox{-\height/2}{\begin{tikzpicture}[scale=.6]
        \node (dots) at (-.5,.3) {$\smash{\vdots}$};
        \node (0) at (-.5,0) {};
        \node (1l) at (1,-1) {};
        \node (1m) at (-.5,-1) {};
        \node (1r) at (-2,-1) {};
        \node (2l) at (-3,-2) {};
        \node (2r) at (-2,-2) {};
        \node (2L) at (-1,-2) {};
        \node (3L) at (-2+.75,-3) {};
        \node (3R) at (1.25-2,-3) {};
        \node (3r) at (.25-2,-3) {};
        \node (3m) at (-2,-3) {};
        \node (3l) at (-2.25,-3) {};
        \draw[color=blue] (1l.center)--(0.center)--(1r.center);
        \draw[color=blue] (1m.center)--(0.center);
        \draw[color=red] (2l.center)--(1r.center);
        \draw[color=red] (1r.center)--(2r.center);
        \draw[color=red] (2L.center)--(1r.center);
        \draw (3l.center)--(2r.center);
        \draw (3m.center)--(2r.center);
        \draw (2r.center)--(3r.center);
        \draw (3L.center)--(2L.center);
        \draw (2L.center)--(3R.center);
        \node (Dots) at (3.5,.3) {$\smash{\vdots}$};
        \node (0') at (3.5,0) {};
        \node (1r') at (2,-1) {};
        \node (1m') at (3.5,-1) {};
        \node (1l') at (5,-1) {};
        \node (2l') at (1,-2) {};
        \node (2r') at (2,-2) {};
        \node (2L') at (3,-2) {};
        \node (3L') at (2.75,-3) {};
        \node (3M') at (3,-3) {};
        \node (3R') at (3.25,-3) {};
        \node (3l') at (2.25,-3) {};
        \node (3r') at (1.75,-3) {};
        \draw[color=blue] (1l'.center)--(0'.center)--(1r'.center);
        \draw[color=blue] (1m'.center)--(0'.center);
        \draw[color=red] (2l'.center)--(1r'.center);
        \draw[color=red] (1r'.center)--(2r'.center);
        \draw[color=red] (2L'.center)--(1r'.center);
        \draw (3l'.center)--(2r'.center);
        \draw (3M'.center)--(2L'.center);
        \draw (2r'.center)--(3r'.center);
        \draw (3L'.center)--(2L'.center);
        \draw (2L'.center)--(3R'.center);
    \end{tikzpicture}}\right)\cdot\left(\raisebox{-\height/2}{\begin{tikzpicture}[scale=.6]
        \node (dots) at (0,.3) {$\smash{\vdots}$};
        \node (0) at (0,0) {};
        \node (1l) at (-.5,-1) {};
        \node (1r) at (.5,-1) {};
        \node (2l) at (-.75,-2) {};
        \node (2m) at (-.5,-2) {};
        \node (2r) at (-.25,-2) {};
        \node (2L) at (.25,-2) {};
        \node (2R) at (.75,-2) {};
        \draw[color=red] (1l.center)--(0.center)--(1r.center);
        \draw (2l.center)--(1l.center);
        \draw (1l.center)--(2r.center);
        \draw (2m.center)--(1l.center);
        \draw (2L.center)--(1r.center);
        \draw (1r.center)--(2R.center);
        \node (dots) at (2,.3) {$\smash{\vdots}$};
        \node (0') at (2,0) {};
        \node (1l') at (1.5,-1) {};
        \node (1r') at (2.5,-1) {};
        \node (2l') at (1.25,-2) {};
        \node (2r') at (1.75,-2) {};
        \node (2L') at (2.25,-2) {};
        \node (2M') at (2.5,-2) {};
        \node (2R') at (2.75,-2) {};
        \draw[color=red] (1l'.center)--(0'.center)--(1r'.center);
        \draw (2l'.center)--(1l'.center);
        \draw (1l'.center)--(2r'.center);
        \draw (2L'.center)--(1r'.center);
        \draw (2M'.center)--(1r'.center);
        \draw (1r'.center)--(2R'.center);
    \end{tikzpicture}}\right)$$
    $$\raisebox{-\height/2}{\begin{tikzpicture}[scale=.6]
        \node (dots) at (0,.3) {$\smash{\vdots}$};
        \node (0) at (0,0) {};
        \node (1l) at (-1,-1) {};
        \node (1r) at (1,-1) {};
        \node (2l) at (2,-2) {};
        \node (2r) at (0,-2) {};
        \node (2L) at (1,-2) {};
        \node (3L) at (.75,-3) {};
        \node (3R) at (1.25,-3) {};
        \node (3r) at (.25,-3) {};
        \node (3m) at (0,-3) {};
        \node (3l) at (-.25,-3) {};
        \draw[color=blue] (1l.center)--(0.center)--(1r.center);
        \draw[color=red] (2l.center)--(1r.center);
        \draw[color=red] (1r.center)--(2r.center);
        \draw[color=red] (2L.center)--(1r.center);
        \draw (3l.center)--(2r.center);
        \draw (3m.center)--(2r.center);
        \draw (2r.center)--(3r.center);
        \draw (3L.center)--(2L.center);
        \draw (2L.center)--(3R.center);
        \node (Dots) at (4,.3) {$\smash{\vdots}$};
        \node (0') at (4,0) {};
        \node (1l') at (3,-1) {};
        \node (1r') at (5,-1) {};
        \node (2l') at (6,-2) {};
        \node (2r') at (4,-2) {};
        \node (2L') at (5,-2) {};
        \node (3L') at (4.75,-3) {};
        \node (3M') at (5,-3) {};
        \node (3R') at (5.25,-3) {};
        \node (3l') at (4.25,-3) {};
        \node (3r') at (3.75,-3) {};
        \draw[color=blue] (1l'.center)--(0'.center)--(1r'.center);
        \draw[color=red] (2l'.center)--(1r'.center);
        \draw[color=red] (1r'.center)--(2r'.center);
        \draw[color=red] (2L'.center)--(1r'.center);
        \draw (3l'.center)--(2r'.center);
        \draw (3M'.center)--(2L'.center);
        \draw (2r'.center)--(3r'.center);
        \draw (3L'.center)--(2L'.center);
        \draw (2L'.center)--(3R'.center);
    \end{tikzpicture}}=\left(\raisebox{-\height/2}{\begin{tikzpicture}[scale=.6]
        \node (dots) at (0,.3) {$\smash{\vdots}$};
        \node (0) at (0,0) {};
        \node (1l) at (-.5,-1) {};
        \node (1r) at (.5,-1) {};
        \node (2l) at (-.75,-2) {};
        \node (2m) at (-.5,-2) {};
        \node (2r) at (-.25,-2) {};
        \node (2L) at (.25,-2) {};
        \node (2R) at (.75,-2) {};
        \draw[color=red] (1l.center)--(0.center)--(1r.center);
        \draw (2l.center)--(1l.center);
        \draw (1l.center)--(2r.center);
        \draw (2m.center)--(1l.center);
        \draw (2L.center)--(1r.center);
        \draw (1r.center)--(2R.center);
        \node (Dots) at (2,.3) {$\smash{\vdots}$};
        \node (0') at (2,0) {};
        \node (1l') at (1.5,-1) {};
        \node (1r') at (2.5,-1) {};
        \node (2l') at (1.25,-2) {};
        \node (2r') at (1.75,-2) {};
        \node (2L') at (2.25,-2) {};
        \node (2M') at (2.5,-2) {};
        \node (2R') at (2.75,-2) {};
        \draw[color=red] (1l'.center)--(0'.center)--(1r'.center);
        \draw (2l'.center)--(1l'.center);
        \draw (1l'.center)--(2r'.center);
        \draw (2L'.center)--(1r'.center);
        \draw (2M'.center)--(1r'.center);
        \draw (1r'.center)--(2R'.center);
    \end{tikzpicture}}\right)\cdot\left(\raisebox{-\height/2}{\begin{tikzpicture}[scale=.6]
        \node (dots) at (-.5,.3) {$\smash{\vdots}$};
        \node (0) at (-.5,0) {};
        \node (1l) at (-2,-1) {};
        \node (1m) at (-.5,-1) {};
        \node (1r) at (1,-1) {};
        \node (2l) at (2,-2) {};
        \node (2r) at (0,-2) {};
        \node (2L) at (1,-2) {};
        \node (3L) at (.75,-3) {};
        \node (3R) at (1.25,-3) {};
        \node (3r) at (.25,-3) {};
        \node (3m) at (0,-3) {};
        \node (3l) at (-.25,-3) {};
        \draw[color=blue] (1l.center)--(0.center)--(1r.center);
        \draw[color=blue] (1m.center)--(0.center);
        \draw[color=red] (2l.center)--(1r.center);
        \draw[color=red] (1r.center)--(2r.center);
        \draw[color=red] (2L.center)--(1r.center);
        \draw (3l.center)--(2r.center);
        \draw (3m.center)--(2r.center);
        \draw (2r.center)--(3r.center);
        \draw (3L.center)--(2L.center);
        \draw (2L.center)--(3R.center);
        \node (Dots) at (3.5,.3) {$\smash{\vdots}$};
        \node (0') at (3.5,0) {};
        \node (1l') at (2,-1) {};
        \node (1m') at (3.5,-1) {};
        \node (1r') at (5,-1) {};
        \node (2l') at (6,-2) {};
        \node (2r') at (4,-2) {};
        \node (2L') at (5,-2) {};
        \node (3L') at (4.75,-3) {};
        \node (3M') at (5,-3) {};
        \node (3R') at (5.25,-3) {};
        \node (3l') at (4.25,-3) {};
        \node (3r') at (3.75,-3) {};
        \draw[color=blue] (1l'.center)--(0'.center)--(1r'.center);
        \draw[color=blue] (1m'.center)--(0'.center);
        \draw[color=red] (2l'.center)--(1r'.center);
        \draw[color=red] (1r'.center)--(2r'.center);
        \draw[color=red] (2L'.center)--(1r'.center);
        \draw (3l'.center)--(2r'.center);
        \draw (3M'.center)--(2L'.center);
        \draw (2r'.center)--(3r'.center);
        \draw (3L'.center)--(2L'.center);
        \draw (2L'.center)--(3R'.center);
    \end{tikzpicture}}\right)\cdot\left(\raisebox{-\height/2}{\begin{tikzpicture}[scale=.6]
        \node (dots) at (0,.3) {$\smash{\vdots}$};
        \node (0) at (0,0) {};
        \node (1l) at (-.5,-1) {};
        \node (1r) at (.5,-1) {};
        \node (2l) at (-.75,-2) {};
        \node (2m) at (.5,-2) {};
        \node (2r) at (-.25,-2) {};
        \node (2L) at (.25,-2) {};
        \node (2R) at (.75,-2) {};
        \draw[color=red] (1l.center)--(0.center)--(1r.center);
        \draw (2l.center)--(1l.center);
        \draw (1l.center)--(2r.center);
        \draw (2m.center)--(1r.center);
        \draw (2L.center)--(1r.center);
        \draw (1r.center)--(2R.center);
        \node (Dots) at (2,.3) {$\smash{\vdots}$};
        \node (0') at (2,0) {};
        \node (1l') at (1.5,-1) {};
        \node (1r') at (2.5,-1) {};
        \node (2l') at (1.25,-2) {};
        \node (2r') at (1.75,-2) {};
        \node (2L') at (2.25,-2) {};
        \node (2M') at (1.5,-2) {};
        \node (2R') at (2.75,-2) {};
        \draw[color=red] (1l'.center)--(0'.center)--(1r'.center);
        \draw (2l'.center)--(1l'.center);
        \draw (1l'.center)--(2r'.center);
        \draw (2L'.center)--(1r'.center);
        \draw (2M'.center)--(1l'.center);
        \draw (1r'.center)--(2R'.center);
    \end{tikzpicture}}\right)$$
    \caption{How to get rid of a $2$-caret just above the hip.}
    \end{figure}
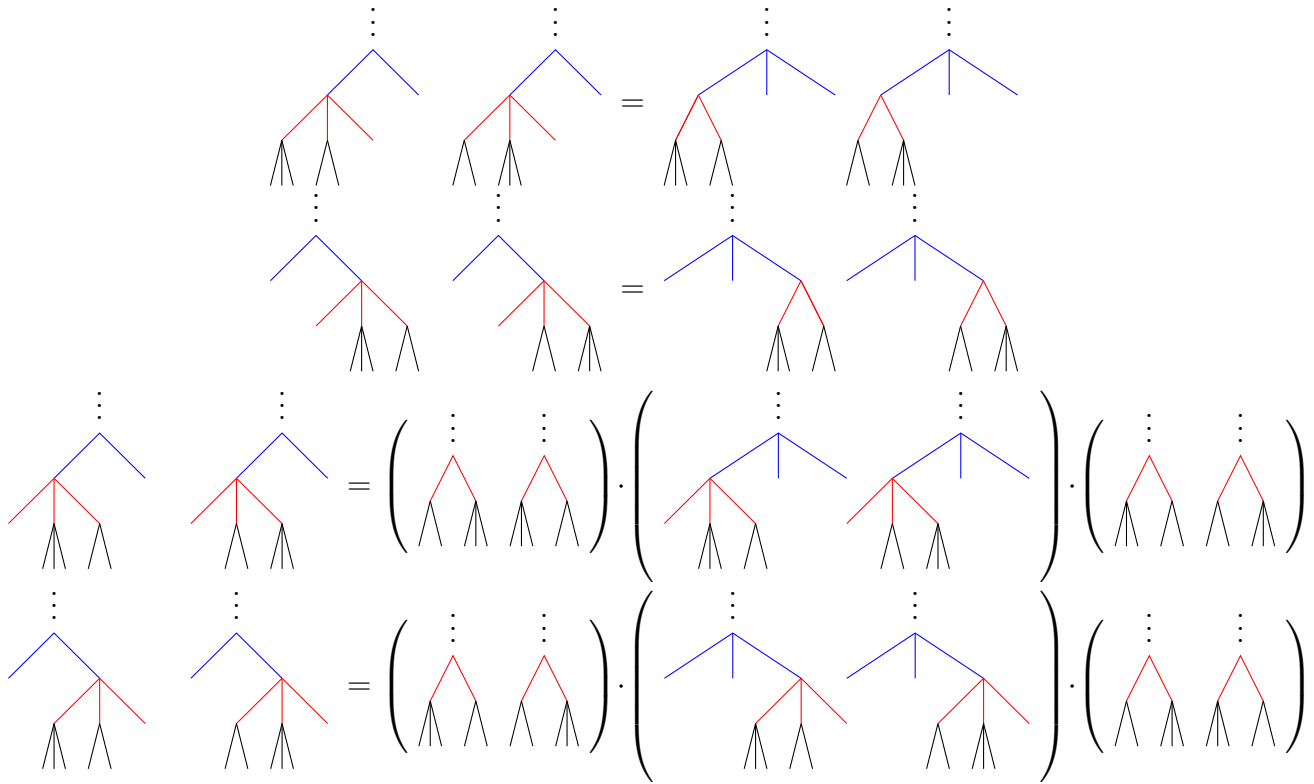

    Each of the factors in these decompositions have one $2$-caret less in their torso, so we can transform them inductively into a product of centipedes, proving that centipedes generate the serpents.
    \end{prf}
    \section{Stick bugs}
    We now introduce the last infinite generating family:

    \begin{dfn} An element is called a left stick bug if it is a centipede where each caret in the torso has the next caret hanging from its left leg, and its feet hang from the two leftmost stems of the hip. A right stick bug is defined analogously. An element is called a stick bug if it is either a left stick bug or a right stick bug.
    \end{dfn}

\begin{figure}[H]
    $$\raisebox{-\height/2}{\begin{tikzpicture}[scale=.6]
        \node (l) at (-4,0) {$l_n\text{:}$};
        \node (0) at (0,0) {};
        \node (1l) at (-1,-1) {};
        \node (1m) at (0,-1) {};
        \node (1r) at (1,-1) {};
        \node (2l) at (-2,-2) {};
        \node (2m) at (-1,-2) {};
        \node (2r) at (0,-2) {};
        \node (nl) at (-3,-3) {};
        \node (ml) at (-4,-4) {};
        \node (mm) at (-3,-4) {};
        \node (mr) at (-2,-4) {};
        \node (fl) at (-4.25,-5) {};
        \node (fm) at (-4,-5) {};
        \node (fr) at (-3.75,-5) {};
        \node (Fl) at (-3.25,-5) {};
        \node (Fr) at (-2.75,-5) {};
        \node (.) at (-2.5,-2.3) {$\iddots$};
        \draw (2l.center)--(1l.center)--(0.center)--(1r.center);
        \draw (1m.center)--(0.center);
        \draw (2r.center)--(1l.center);
        \draw (2m.center)--(1l.center);
        \draw (ml.center)--(nl.center);
        \draw (mm.center)--(nl.center);
        \draw (mr.center)--(nl.center);
        \draw (fl.center)--(ml.center);
        \draw (fm.center)--(ml.center);
        \draw (fr.center)--(ml.center);
        \draw (Fl.center)--(mm.center);
        \draw (Fr.center)--(mm.center);
        \draw[decoration={brace,mirror,raise=5pt},decorate] (0.center)--node[above left=4pt and 4pt] {$n$} (nl.center);
        \node (0') at (4,0) {};
        \node (1l') at (3,-1) {};
        \node (1m') at (4,-1) {};
        \node (1r') at (5,-1) {};
        \node (2l') at (2,-2) {};
        \node (2m') at (3,-2) {};
        \node (2r') at (4,-2) {};
        \node (nl') at (1,-3) {};
        \node (ml') at (0,-4) {};
        \node (mm') at (1,-4) {};
        \node (mr') at (2,-4) {};
        \node (fl') at (-.25,-5) {};
        \node (fr') at (.25,-5) {};
        \node (Fl') at (.75,-5) {};
        \node (Fm') at (1,-5) {};
        \node (Fr') at (1.25,-5) {};
        \node (.') at (1.5,-2.3) {$\iddots$};
        \draw (2l'.center)--(1l'.center)--(0'.center)--(1r'.center);
        \draw (1m'.center)--(0'.center);
        \draw (2r'.center)--(1l'.center);
        \draw (2m'.center)--(1l'.center);
        \draw (ml'.center)--(nl'.center);
        \draw (mm'.center)--(nl'.center);
        \draw (mr'.center)--(nl'.center);
        \draw (fl'.center)--(ml'.center);
        \draw (fr'.center)--(ml'.center);
        \draw (Fl'.center)--(mm'.center);
        \draw (Fm'.center)--(mm'.center);
        \draw (Fr'.center)--(mm'.center);
        \draw[decoration={brace,mirror,raise=5pt},decorate] (0'.center)--node[above left=4pt and 4pt] {$n$} (nl'.center);
    \end{tikzpicture}}\qquad\raisebox{-\height/2}{\begin{tikzpicture}[scale=.6]
        \node (l) at (-6,0) {$r_n\text{:}$};
        \node (0) at (0,0) {};
        \node (1l) at (1,-1) {};
        \node (1m) at (0,-1) {};
        \node (1r) at (-1,-1) {};
        \node (2l) at (2,-2) {};
        \node (2m) at (1,-2) {};
        \node (2r) at (0,-2) {};
        \node (nl) at (3,-3) {};
        \node (ml) at (4,-4) {};
        \node (mm) at (3,-4) {};
        \node (mr) at (2,-4) {};
        \node (fl) at (4.25,-5) {};
        \node (fm) at (4,-5) {};
        \node (fr) at (3.75,-5) {};
        \node (Fl) at (3.25,-5) {};
        \node (Fr) at (2.75,-5) {};
        \node (.) at (2.5,-2.3) {$\ddots$};
        \draw (2l.center)--(1l.center)--(0.center)--(1r.center);
        \draw (1m.center)--(0.center);
        \draw (2r.center)--(1l.center);
        \draw (2m.center)--(1l.center);
        \draw (ml.center)--(nl.center);
        \draw (mm.center)--(nl.center);
        \draw (mr.center)--(nl.center);
        \draw (fl.center)--(ml.center);
        \draw (fm.center)--(ml.center);
        \draw (fr.center)--(ml.center);
        \draw (Fl.center)--(mm.center);
        \draw (Fr.center)--(mm.center);
        \draw[decoration={brace,raise=5pt},decorate] (0.center)--node[above right=4pt and 4pt] {$n$} (nl.center);
        \node (0') at (-4,0) {};
        \node (1l') at (-3,-1) {};
        \node (1m') at (-4,-1) {};
        \node (1r') at (-5,-1) {};
        \node (2l') at (-2,-2) {};
        \node (2m') at (-3,-2) {};
        \node (2r') at (-4,-2) {};
        \node (nl') at (-1,-3) {};
        \node (ml') at (0,-4) {};
        \node (mm') at (-1,-4) {};
        \node (mr') at (-2,-4) {};
        \node (fl') at (.25,-5) {};
        \node (fr') at (-.25,-5) {};
        \node (Fl') at (-.75,-5) {};
        \node (Fm') at (-1,-5) {};
        \node (Fr') at (-1.25,-5) {};
        \node (.') at (-1.5,-2.3) {$\ddots$};
        \draw (2l'.center)--(1l'.center)--(0'.center)--(1r'.center);
        \draw (1m'.center)--(0'.center);
        \draw (2r'.center)--(1l'.center);
        \draw (2m'.center)--(1l'.center);
        \draw (ml'.center)--(nl'.center);
        \draw (mm'.center)--(nl'.center);
        \draw (mr'.center)--(nl'.center);
        \draw (fl'.center)--(ml'.center);
        \draw (fr'.center)--(ml'.center);
        \draw (Fl'.center)--(mm'.center);
        \draw (Fm'.center)--(mm'.center);
        \draw (Fr'.center)--(mm'.center);
        \draw[decoration={brace,raise=5pt},decorate] (0'.center)--node[above right=4pt and 4pt] {$n$} (nl'.center);
    \end{tikzpicture}}$$
    \centering
    \includegraphics[scale=0.5]{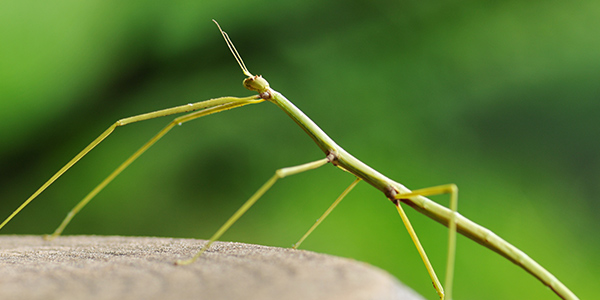}
    \caption{The stick bugs $l_n$ and $r_n$, and a photo of a real stick bug for comparison. Source: https://www.nwf.org/Educational-Resources/Wildlife-Guide/Invertebrates/Walking-Sticks}
    \label{getstickbuggedlol}
    \end{figure}

    There are only four stick bugs with a given number $n$ of torso pieces. Two of them (see \Cref{getstickbuggedlol}) are given the names $l_n$ and $r_n$, and the other two are their inverses, which have their foot carets the other way around.
    The stick bugs $l_0$ and $r_0$, according to \Cref{fg}, are going to be the two generators of the group. For convenience, we may drop the subindex and simply call them $l$ and $r$.\\

    Our goal now is to prove the following:

    \begin{prp}\label{stickgen}
        The stick bugs generate $F\left(\frac32\right)$.
    \end{prp}

    And we start with a lemma:

    \begin{lmm}\label{xiv}
        The relations in \Cref{xiv2} hold.
    \end{lmm}
    \begin{prf}
        It is a lengthy albeit straight-forward calculation (for a computer).
    \end{prf}

    \begin{figure}[H]
        \begin{align*}\raisebox{-\height/2}{\begin{tikzpicture}[scale=.6]
        \node (0) at (0,0) {};
        \node (1l) at (-1,-1) {};
        \node (1m) at (0,-1) {};
        \node (1r) at (1,-1) {};
        \node (2l) at (-1,-2) {};
        \node (2m) at (0,-2) {};
        \node (2r) at (1,-2) {};
        \node (fl) at (-1.25,-3) {};
        \node (fm) at (-1,-3) {};
        \node (fr) at (-0.75,-3) {};
        \node (Fl) at (0.25,-3) {};
        \node (Fr) at (-0.25,-3) {};
        \draw (1l.center)--(0.center)--(1r.center);
        \draw (1m.center)--(0.center);
        \draw (2r.center)--(1m.center);
        \draw (2m.center)--(1m.center);
        \draw (2l.center)--(1m.center);
        \draw (fl.center)--(2l.center);
        \draw (fm.center)--(2l.center);
        \draw (fr.center)--(2l.center);
        \draw (Fl.center)--(2m.center);
        \draw (Fr.center)--(2m.center);
        \node (0') at (3,0) {};
        \node (1l') at (2,-1) {};
        \node (1m') at (3,-1) {};
        \node (1r') at (4,-1) {};
        \node (2l') at (2,-2) {};
        \node (2m') at (3,-2) {};
        \node (2r') at (4,-2) {};
        \node (fl') at (1.75,-3) {};
        \node (fr') at (2.25,-3) {};
        \node (Fl') at (3.25,-3) {};
        \node (Fm') at (3,-3) {};
        \node (Fr') at (2.75,-3) {};
        \draw (1l'.center)--(0'.center)--(1r'.center);
        \draw (1m'.center)--(0'.center);
        \draw (2r'.center)--(1m'.center);
        \draw (2m'.center)--(1m'.center);
        \draw (2l'.center)--(1m'.center);
        \draw (fl'.center)--(2l'.center);
        \draw (fr'.center)--(2l'.center);
        \draw (Fl'.center)--(2m'.center);
        \draw (Fm'.center)--(2m'.center);
        \draw (Fr'.center)--(2m'.center);
    \end{tikzpicture}}&=r^{-1}l^{-1}r^{-1}lrl^{-1}r^{-1}lrl^{-1}r^{-1}l^{-1}rlr^{-1}l^{-1}rlrl=\left[[r,l]^2\cdot r\large\text,l\right]\\
    \raisebox{-\height/2}{\begin{tikzpicture}[scale=.6]
        \node (0) at (0,0) {};
        \node (1l) at (1,-1) {};
        \node (1m) at (0,-1) {};
        \node (1r) at (-1,-1) {};
        \node (2l) at (1,-2) {};
        \node (2m) at (0,-2) {};
        \node (2r) at (-1,-2) {};
        \node (fl) at (1.25,-3) {};
        \node (fr) at (0.75,-3) {};
        \node (Fl) at (-0.25,-3) {};
        \node (Fm) at (0,-3) {};
        \node (Fr) at (0.25,-3) {};
        \draw (1l.center)--(0.center)--(1r.center);
        \draw (1m.center)--(0.center);
        \draw (2r.center)--(1m.center);
        \draw (2m.center)--(1m.center);
        \draw (2l.center)--(1m.center);
        \draw (fl.center)--(2l.center);
        \draw (fr.center)--(2l.center);
        \draw (Fl.center)--(2m.center);
        \draw (Fm.center)--(2m.center);
        \draw (Fr.center)--(2m.center);
        \node (0') at (3,0) {};
        \node (1l') at (2,-1) {};
        \node (1m') at (3,-1) {};
        \node (1r') at (4,-1) {};
        \node (2l') at (2,-2) {};
        \node (2m') at (3,-2) {};
        \node (2r') at (4,-2) {};
        \node (fl') at (2.75,-3) {};
        \node (fr') at (3.25,-3) {};
        \node (Fl') at (4.25,-3) {};
        \node (Fm') at (4,-3) {};
        \node (Fr') at (3.75,-3) {};
        \draw (1l'.center)--(0'.center)--(1r'.center);
        \draw (1m'.center)--(0'.center);
        \draw (2r'.center)--(1m'.center);
        \draw (2m'.center)--(1m'.center);
        \draw (2l'.center)--(1m'.center);
        \draw (fl'.center)--(2m'.center);
        \draw (fr'.center)--(2m'.center);
        \draw (Fl'.center)--(2r'.center);
        \draw (Fm'.center)--(2r'.center);
        \draw (Fr'.center)--(2r'.center);
    \end{tikzpicture}}&=rlrlr^{-1}l^{-1}rlr^{-1}l^{-1}r^{-1}lrl^{-1}r^{-1}lrl^{-1}r^{-1}l^{-1}=\left[r^{-1}{\large\text,}[l^{-1},r^{-1}]^2\cdot l^{-1}\right]\end{align*}
    \caption{Writing two particular elements using $l$ and $r$.}
    \label{xiv2}
    \end{figure}

    The fact that these two elements are in the span of $l$ and $r$ will be important now as well as later, and was first indirectly proven by the authors with computer assistance. Then, explicit longer words in $l$ and $r$ were found and then replaced by the current $20$-long words, also with the aid of a computer. Shorter words for these elements do not exist.

    Note that by adding any torso on top of every tree in the above product we obtain the following corollary:

    \begin{crl}\label{indstick} The centipedes whose hip hangs from the middle stem of a $3$-caret can be written in terms of other centipedes of strictly less depth.
    \end{crl}

    By nesting the relations with themselves, we get another corollary that will be useful later.

    \begin{crl}\label{xvnested} Any centipede in which every non-foot caret hangs from the middle stem of the previous one, belongs to the span of $l$ and $r$.
    \end{crl}

    We can now prove that the stick bugs generate the group. We will do that by generating every centipede:\\

    \begin{prf}[ of \Cref{stickgen}]
    We need to write a centipede as a product of stick bugs and we will use induction on the depth, or equivalently, the number of torso pieces. Notice that the only centipedes without torso are the stick bugs $l$ and $r$ and their inverses, so we can assume that there is at least a torso piece and that hip is hanging. If the hip hangs from the middle, \Cref{indstick} together with the induction hypothesis allows us to write the element as a product of stick bugs.

    Now suppose that, even though the hip does not hang from a middle stem, there is a caret in the torso that does. In that case, via a series of conjugations by shorter centipedes (see \Cref{redreçantcentpeus}) we are able to make the following caret hang from the middle stem as well. By repeating this, we are able to make a chain reaction that ends with the hip hanging from the middle stem, reducing it to the previous case.

    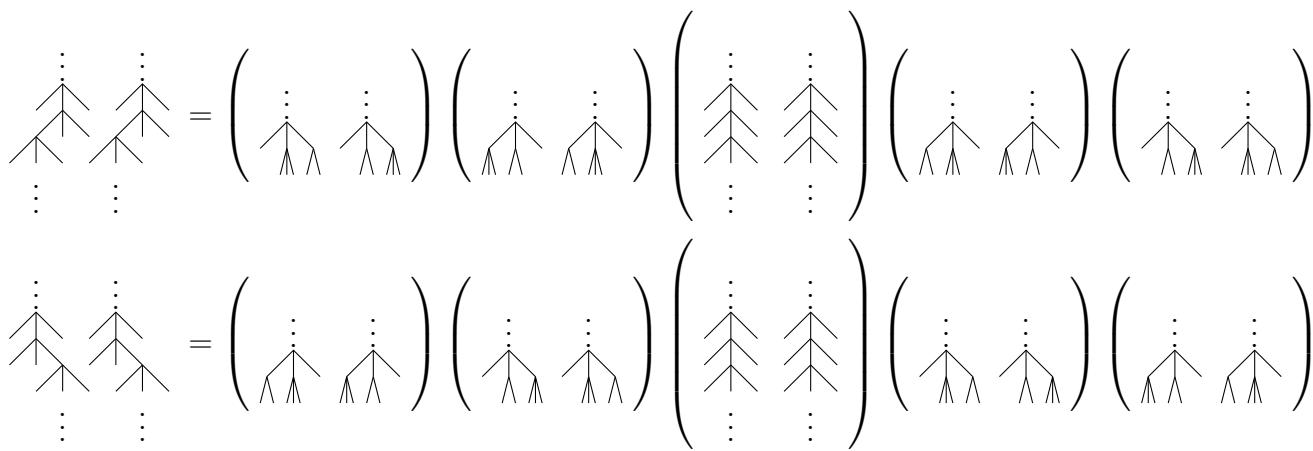
\begin{figure}[H]
    $$\raisebox{-\height/2}{\begin{tikzpicture}[scale=.35]
        \node (-1) at (0,1) {$\vdots$};
        \node (0) at (0,0) {};
        \node (1l) at (-1,-1) {};
        \node (1m) at (0,-1) {};
        \node (1r) at (1,-1) {};
        \node (2l) at (-1,-2) {};
        \node (2m) at (0,-2) {};
        \node (2r) at (1,-2) {};
        \node (3l) at (-2,-3) {};
        \node (3m) at (-1,-3) {};
        \node (3r) at (0,-3) {};
        \node (4m) at (-1,-4) {$\vdots$};
        \draw (1l.center)--(0.center)--(1r.center);
        \draw (1m.center)--(0.center);
        \draw (2r.center)--(1m.center);
        \draw (2m.center)--(1m.center);
        \draw (2l.center)--(1m.center);
        \draw (3l.center)--(2l.center);
        \draw (3m.center)--(2l.center);
        \draw (3r.center)--(2l.center);
        \node (-1') at (3,1) {$\vdots$};
        \node (0') at (3,0) {};
        \node (1l') at (2,-1) {};
        \node (1m') at (3,-1) {};
        \node (1r') at (4,-1) {};
        \node (2l') at (2,-2) {};
        \node (2m') at (3,-2) {};
        \node (2r') at (4,-2) {};
        \node (2L') at (2,-2) {};
        \node (2R') at (4,-2) {};
        \node (3l') at (1,-3) {};
        \node (3m') at (2,-3) {};
        \node (3r') at (3,-3) {};
        \node (4m') at (2,-4) {$\vdots$};
        \draw (1l'.center)--(0'.center)--(1r'.center);
        \draw (1m'.center)--(0'.center);
        \draw (2r'.center)--(1m'.center);
        \draw (2m'.center)--(1m'.center);
        \draw (2l'.center)--(1m'.center);
        \draw (3l'.center)--(2l'.center);
        \draw (3m'.center)--(2l'.center);
        \draw (3r'.center)--(2l'.center);
    \end{tikzpicture}}=\left(\raisebox{-\height/2}{\begin{tikzpicture}[scale=.35]
        \node (-1) at (0,1) {$\vdots$};
        \node (0) at (0,0) {};
        \node (1l) at (-1,-1) {};
        \node (1m) at (0,-1) {};
        \node (1r) at (1,-1) {};
        \node (2l) at (-.25,-2) {};
        \node (2m) at (0,-2) {};
        \node (2r) at (.25,-2) {};
        \node (2L) at (.75,-2) {};
        \node (2R) at (1.25,-2) {};
        \draw (1l.center)--(0.center)--(1r.center);
        \draw (1m.center)--(0.center);
        \draw (2r.center)--(1m.center);
        \draw (2m.center)--(1m.center);
        \draw (2l.center)--(1m.center);
        \draw (2L.center)--(1r.center);
        \draw (2R.center)--(1r.center);
        \node (-1') at (3,1) {$\vdots$};
        \node (0') at (3,0) {};
        \node (1l') at (2,-1) {};
        \node (1m') at (3,-1) {};
        \node (1r') at (4,-1) {};
        \node (2l') at (2.75,-2) {};
        \node (2r') at (3.25,-2) {};
        \node (2L') at (3.75,-2) {};
        \node (2M') at (4,-2) {};
        \node (2R') at (4.25,-2) {};
        \draw (1l'.center)--(0'.center)--(1r'.center);
        \draw (1m'.center)--(0'.center);
        \draw (2r'.center)--(1m'.center);
        \draw (2l'.center)--(1m'.center);
        \draw (2L'.center)--(1r'.center);
        \draw (2M'.center)--(1r'.center);
        \draw (2R'.center)--(1r'.center);
    \end{tikzpicture}}\right)
    \left(\raisebox{-\height/2}{\begin{tikzpicture}[scale=.35]
        \node (-1) at (0,1) {$\vdots$};
        \node (0) at (0,0) {};
        \node (1l) at (-1,-1) {};
        \node (1m) at (0,-1) {};
        \node (1r) at (1,-1) {};
        \node (2l) at (-.25,-2) {};
        \node (2r) at (.25,-2) {};
        \node (2L) at (-1.25,-2) {};
        \node (2M) at (-1,-2) {};
        \node (2R) at (-.75,-2) {};
        \draw (1l.center)--(0.center)--(1r.center);
        \draw (1m.center)--(0.center);
        \draw (2r.center)--(1m.center);
        \draw (2l.center)--(1m.center);
        \draw (2L.center)--(1l.center);
        \draw (2M.center)--(1l.center);
        \draw (2R.center)--(1l.center);
        \node (-1') at (3,1) {$\vdots$};
        \node (0') at (3,0) {};
        \node (1l') at (2,-1) {};
        \node (1m') at (3,-1) {};
        \node (1r') at (4,-1) {};
        \node (2l') at (2.75,-2) {};
        \node (2m') at (3,-2) {};
        \node (2r') at (3.25,-2) {};
        \node (2L') at (1.75,-2) {};
        \node (2R') at (2.25,-2) {};
        \draw (1l'.center)--(0'.center)--(1r'.center);
        \draw (1m'.center)--(0'.center);
        \draw (2r'.center)--(1m'.center);
        \draw (2m'.center)--(1m'.center);
        \draw (2l'.center)--(1m'.center);
        \draw (2L'.center)--(1l'.center);
        \draw (2R'.center)--(1l'.center);
    \end{tikzpicture}}\right)
    \left(\raisebox{-\height/2}{\begin{tikzpicture}[scale=.35]
        \node (-1) at (0,1) {$\vdots$};
        \node (0) at (0,0) {};
        \node (1l) at (-1,-1) {};
        \node (1m) at (0,-1) {};
        \node (1r) at (1,-1) {};
        \node (2l) at (-1,-2) {};
        \node (2m) at (0,-2) {};
        \node (2r) at (1,-2) {};
        \node (3l) at (-1,-3) {};
        \node (3m) at (0,-3) {};
        \node (3r) at (1,-3) {};
        \node (4m) at (0,-4) {$\vdots$};
        \draw (1l.center)--(0.center)--(1r.center);
        \draw (1m.center)--(0.center);
        \draw (2r.center)--(1m.center);
        \draw (2m.center)--(1m.center);
        \draw (2l.center)--(1m.center);
        \draw (3l.center)--(2m.center);
        \draw (3m.center)--(2m.center);
        \draw (3r.center)--(2m.center);
        \node (-1') at (3,1) {$\vdots$};
        \node (0') at (3,0) {};
        \node (1l') at (2,-1) {};
        \node (1m') at (3,-1) {};
        \node (1r') at (4,-1) {};
        \node (2l') at (2,-2) {};
        \node (2m') at (3,-2) {};
        \node (2r') at (4,-2) {};
        \node (3l') at (2,-3) {};
        \node (3m') at (3,-3) {};
        \node (3r') at (4,-3) {};
        \node (4m') at (3,-4) {$\vdots$};
        \draw (1l'.center)--(0'.center)--(1r'.center);
        \draw (1m'.center)--(0'.center);
        \draw (2r'.center)--(1m'.center);
        \draw (2m'.center)--(1m'.center);
        \draw (2l'.center)--(1m'.center);
        \draw (3l'.center)--(2m'.center);
        \draw (3m'.center)--(2m'.center);
        \draw (3r'.center)--(2m'.center);
    \end{tikzpicture}}\right)\left(\raisebox{-\height/2}{\begin{tikzpicture}[scale=.35]
        \node (-1) at (0,1) {$\vdots$};
        \node (0) at (0,0) {};
        \node (1l) at (-1,-1) {};
        \node (1m) at (0,-1) {};
        \node (1r) at (1,-1) {};
        \node (2l) at (-.25,-2) {};
        \node (2m) at (0,-2) {};
        \node (2r) at (.25,-2) {};
        \node (2L) at (-1.25,-2) {};
        \node (2R) at (-.75,-2) {};
        \draw (1l.center)--(0.center)--(1r.center);
        \draw (1m.center)--(0.center);
        \draw (2r.center)--(1m.center);
        \draw (2m.center)--(1m.center);
        \draw (2l.center)--(1m.center);
        \draw (2L.center)--(1l.center);
        \draw (2R.center)--(1l.center);
        \node (-1') at (3,1) {$\vdots$};
        \node (0') at (3,0) {};
        \node (1l') at (2,-1) {};
        \node (1m') at (3,-1) {};
        \node (1r') at (4,-1) {};
        \node (2l') at (2.75,-2) {};
        \node (2r') at (3.25,-2) {};
        \node (2L') at (1.75,-2) {};
        \node (2M') at (2,-2) {};
        \node (2R') at (2.25,-2) {};
        \draw (1l'.center)--(0'.center)--(1r'.center);
        \draw (1m'.center)--(0'.center);
        \draw (2r'.center)--(1m'.center);
        \draw (2l'.center)--(1m'.center);
        \draw (2L'.center)--(1l'.center);
        \draw (2M'.center)--(1l'.center);
        \draw (2R'.center)--(1l'.center);
    \end{tikzpicture}}\right)\left(\raisebox{-\height/2}{\begin{tikzpicture}[scale=.35]
        \node (-1) at (0,1) {$\vdots$};
        \node (0) at (0,0) {};
        \node (1l) at (-1,-1) {};
        \node (1m) at (0,-1) {};
        \node (1r) at (1,-1) {};
        \node (2l) at (-.25,-2) {};
        \node (2r) at (.25,-2) {};
        \node (2L) at (.75,-2) {};
        \node (2M) at (1,-2) {};
        \node (2R) at (1.25,-2) {};
        \draw (1l.center)--(0.center)--(1r.center);
        \draw (1m.center)--(0.center);
        \draw (2r.center)--(1m.center);
        \draw (2l.center)--(1m.center);
        \draw (2L.center)--(1r.center);
        \draw (2M.center)--(1r.center);
        \draw (2R.center)--(1r.center);
        \node (-1') at (3,1) {$\vdots$};
        \node (0') at (3,0) {};
        \node (1l') at (2,-1) {};
        \node (1m') at (3,-1) {};
        \node (1r') at (4,-1) {};
        \node (2l') at (2.75,-2) {};
        \node (2m') at (3,-2) {};
        \node (2r') at (3.25,-2) {};
        \node (2L') at (3.75,-2) {};
        \node (2R') at (4.25,-2) {};
        \draw (1l'.center)--(0'.center)--(1r'.center);
        \draw (1m'.center)--(0'.center);
        \draw (2r'.center)--(1m'.center);
        \draw (2m'.center)--(1m'.center);
        \draw (2l'.center)--(1m'.center);
        \draw (2L'.center)--(1r'.center);
        \draw (2R'.center)--(1r'.center);
    \end{tikzpicture}}\right)$$
    $$\raisebox{-\height/2}{\begin{tikzpicture}[scale=.35]
        \node (-1) at (0,1) {$\vdots$};
        \node (0) at (0,0) {};
        \node (1l) at (1,-1) {};
        \node (1m) at (0,-1) {};
        \node (1r) at (-1,-1) {};
        \node (2l) at (1,-2) {};
        \node (2m) at (0,-2) {};
        \node (2r) at (-1,-2) {};
        \node (3l) at (2,-3) {};
        \node (3m) at (1,-3) {};
        \node (3r) at (0,-3) {};
        \node (4m) at (1,-4) {$\vdots$};
        \draw (1l.center)--(0.center)--(1r.center);
        \draw (1m.center)--(0.center);
        \draw (2r.center)--(1m.center);
        \draw (2m.center)--(1m.center);
        \draw (2l.center)--(1m.center);
        \draw (3l.center)--(2l.center);
        \draw (3m.center)--(2l.center);
        \draw (3r.center)--(2l.center);
        \node (-1') at (-3,1) {$\vdots$};
        \node (0') at (-3,0) {};
        \node (1l') at (-2,-1) {};
        \node (1m') at (-3,-1) {};
        \node (1r') at (-4,-1) {};
        \node (2l') at (-2,-2) {};
        \node (2m') at (-3,-2) {};
        \node (2r') at (-4,-2) {};
        \node (2L') at (-2,-2) {};
        \node (2R') at (-4,-2) {};
        \node (3l') at (-1,-3) {};
        \node (3m') at (-2,-3) {};
        \node (3r') at (-3,-3) {};
        \node (4m') at (-2,-4) {$\vdots$};
        \draw (1l'.center)--(0'.center)--(1r'.center);
        \draw (1m'.center)--(0'.center);
        \draw (2r'.center)--(1m'.center);
        \draw (2m'.center)--(1m'.center);
        \draw (2l'.center)--(1m'.center);
        \draw (3l'.center)--(2l'.center);
        \draw (3m'.center)--(2l'.center);
        \draw (3r'.center)--(2l'.center);
    \end{tikzpicture}}=\left(\raisebox{-\height/2}{\begin{tikzpicture}[scale=.35]
        \node (-1) at (0,1) {$\vdots$};
        \node (0) at (0,0) {};
        \node (1l) at (-1,-1) {};
        \node (1m) at (0,-1) {};
        \node (1r) at (1,-1) {};
        \node (2l) at (-.25,-2) {};
        \node (2m) at (0,-2) {};
        \node (2r) at (.25,-2) {};
        \node (2L) at (-1.25,-2) {};
        \node (2R) at (-.75,-2) {};
        \draw (1l.center)--(0.center)--(1r.center);
        \draw (1m.center)--(0.center);
        \draw (2r.center)--(1m.center);
        \draw (2m.center)--(1m.center);
        \draw (2l.center)--(1m.center);
        \draw (2L.center)--(1l.center);
        \draw (2R.center)--(1l.center);
        \node (-1') at (3,1) {$\vdots$};
        \node (0') at (3,0) {};
        \node (1l') at (2,-1) {};
        \node (1m') at (3,-1) {};
        \node (1r') at (4,-1) {};
        \node (2l') at (2.75,-2) {};
        \node (2r') at (3.25,-2) {};
        \node (2L') at (1.75,-2) {};
        \node (2M') at (2,-2) {};
        \node (2R') at (2.25,-2) {};
        \draw (1l'.center)--(0'.center)--(1r'.center);
        \draw (1m'.center)--(0'.center);
        \draw (2r'.center)--(1m'.center);
        \draw (2l'.center)--(1m'.center);
        \draw (2L'.center)--(1l'.center);
        \draw (2M'.center)--(1l'.center);
        \draw (2R'.center)--(1l'.center);
    \end{tikzpicture}}\right)\left(\raisebox{-\height/2}{\begin{tikzpicture}[scale=.35]
        \node (-1) at (0,1) {$\vdots$};
        \node (0) at (0,0) {};
        \node (1l) at (-1,-1) {};
        \node (1m) at (0,-1) {};
        \node (1r) at (1,-1) {};
        \node (2l) at (-.25,-2) {};
        \node (2r) at (.25,-2) {};
        \node (2L) at (.75,-2) {};
        \node (2M) at (1,-2) {};
        \node (2R) at (1.25,-2) {};
        \draw (1l.center)--(0.center)--(1r.center);
        \draw (1m.center)--(0.center);
        \draw (2r.center)--(1m.center);
        \draw (2l.center)--(1m.center);
        \draw (2L.center)--(1r.center);
        \draw (2M.center)--(1r.center);
        \draw (2R.center)--(1r.center);
        \node (-1') at (3,1) {$\vdots$};
        \node (0') at (3,0) {};
        \node (1l') at (2,-1) {};
        \node (1m') at (3,-1) {};
        \node (1r') at (4,-1) {};
        \node (2l') at (2.75,-2) {};
        \node (2m') at (3,-2) {};
        \node (2r') at (3.25,-2) {};
        \node (2L') at (3.75,-2) {};
        \node (2R') at (4.25,-2) {};
        \draw (1l'.center)--(0'.center)--(1r'.center);
        \draw (1m'.center)--(0'.center);
        \draw (2r'.center)--(1m'.center);
        \draw (2m'.center)--(1m'.center);
        \draw (2l'.center)--(1m'.center);
        \draw (2L'.center)--(1r'.center);
        \draw (2R'.center)--(1r'.center);
    \end{tikzpicture}}\right)\left(\raisebox{-\height/2}{\begin{tikzpicture}[scale=.35]
        \node (-1) at (0,1) {$\vdots$};
        \node (0) at (0,0) {};
        \node (1l) at (-1,-1) {};
        \node (1m) at (0,-1) {};
        \node (1r) at (1,-1) {};
        \node (2l) at (-1,-2) {};
        \node (2m) at (0,-2) {};
        \node (2r) at (1,-2) {};
        \node (3l) at (-1,-3) {};
        \node (3m) at (0,-3) {};
        \node (3r) at (1,-3) {};
        \node (4m) at (0,-4) {$\vdots$};
        \draw (1l.center)--(0.center)--(1r.center);
        \draw (1m.center)--(0.center);
        \draw (2r.center)--(1m.center);
        \draw (2m.center)--(1m.center);
        \draw (2l.center)--(1m.center);
        \draw (3l.center)--(2m.center);
        \draw (3m.center)--(2m.center);
        \draw (3r.center)--(2m.center);
        \node (-1') at (3,1) {$\vdots$};
        \node (0') at (3,0) {};
        \node (1l') at (2,-1) {};
        \node (1m') at (3,-1) {};
        \node (1r') at (4,-1) {};
        \node (2l') at (2,-2) {};
        \node (2m') at (3,-2) {};
        \node (2r') at (4,-2) {};
        \node (2L') at (2,-2) {};
        \node (2R') at (4,-2) {};
        \node (3l') at (2,-3) {};
        \node (3m') at (3,-3) {};
        \node (3r') at (4,-3) {};
        \node (4m') at (3,-4) {$\vdots$};
        \draw (1l'.center)--(0'.center)--(1r'.center);
        \draw (1m'.center)--(0'.center);
        \draw (2r'.center)--(1m'.center);
        \draw (2m'.center)--(1m'.center);
        \draw (2l'.center)--(1m'.center);
        \draw (3l'.center)--(2m'.center);
        \draw (3m'.center)--(2m'.center);
        \draw (3r'.center)--(2m'.center);
    \end{tikzpicture}}\right)\left(\raisebox{-\height/2}{\begin{tikzpicture}[scale=.35]
        \node (-1) at (0,1) {$\vdots$};
        \node (0) at (0,0) {};
        \node (1l) at (-1,-1) {};
        \node (1m) at (0,-1) {};
        \node (1r) at (1,-1) {};
        \node (2l) at (-.25,-2) {};
        \node (2m) at (0,-2) {};
        \node (2r) at (.25,-2) {};
        \node (2L) at (.75,-2) {};
        \node (2R) at (1.25,-2) {};
        \draw (1l.center)--(0.center)--(1r.center);
        \draw (1m.center)--(0.center);
        \draw (2r.center)--(1m.center);
        \draw (2m.center)--(1m.center);
        \draw (2l.center)--(1m.center);
        \draw (2L.center)--(1r.center);
        \draw (2R.center)--(1r.center);
        \node (-1') at (3,1) {$\vdots$};
        \node (0') at (3,0) {};
        \node (1l') at (2,-1) {};
        \node (1m') at (3,-1) {};
        \node (1r') at (4,-1) {};
        \node (2l') at (2.75,-2) {};
        \node (2r') at (3.25,-2) {};
        \node (2L') at (3.75,-2) {};
        \node (2M') at (4,-2) {};
        \node (2R') at (4.25,-2) {};
        \draw (1l'.center)--(0'.center)--(1r'.center);
        \draw (1m'.center)--(0'.center);
        \draw (2r'.center)--(1m'.center);
        \draw (2l'.center)--(1m'.center);
        \draw (2L'.center)--(1r'.center);
        \draw (2M'.center)--(1r'.center);
        \draw (2R'.center)--(1r'.center);
    \end{tikzpicture}}\right)
    \left(\raisebox{-\height/2}{\begin{tikzpicture}[scale=.35]
        \node (-1) at (0,1) {$\vdots$};
        \node (0) at (0,0) {};
        \node (1l) at (-1,-1) {};
        \node (1m) at (0,-1) {};
        \node (1r) at (1,-1) {};
        \node (2l) at (-.25,-2) {};
        \node (2r) at (.25,-2) {};
        \node (2L) at (-1.25,-2) {};
        \node (2M) at (-1,-2) {};
        \node (2R) at (-.75,-2) {};
        \draw (1l.center)--(0.center)--(1r.center);
        \draw (1m.center)--(0.center);
        \draw (2r.center)--(1m.center);
        \draw (2l.center)--(1m.center);
        \draw (2L.center)--(1l.center);
        \draw (2M.center)--(1l.center);
        \draw (2R.center)--(1l.center);
        \node (-1') at (3,1) {$\vdots$};
        \node (0') at (3,0) {};
        \node (1l') at (2,-1) {};
        \node (1m') at (3,-1) {};
        \node (1r') at (4,-1) {};
        \node (2l') at (2.75,-2) {};
        \node (2m') at (3,-2) {};
        \node (2r') at (3.25,-2) {};
        \node (2L') at (1.75,-2) {};
        \node (2R') at (2.25,-2) {};
        \draw (1l'.center)--(0'.center)--(1r'.center);
        \draw (1m'.center)--(0'.center);
        \draw (2r'.center)--(1m'.center);
        \draw (2m'.center)--(1m'.center);
        \draw (2l'.center)--(1m'.center);
        \draw (2L'.center)--(1l'.center);
        \draw (2R'.center)--(1l'.center);
    \end{tikzpicture}}\right)$$
    \caption{Conjugating by centipedes with less torso pieces we can make all subsequent carets (including the hip) hang from the middle stem.}
    \label{redreçantcentpeus}
    \end{figure}

    Hanging from the middle stem is now ruled out, so we study what happens if there are carets hanging from left stems and other carets hanging from right stems. In that case, there will be a place where both situations occur one after the other. Again, via conjugation by shorter centipedes (see \Cref{mesconjs}), we are able to reduce that to the previous case.

    \begin{figure}[H]
    $$\raisebox{-\height/2}{\begin{tikzpicture}[scale=.5]
        \node (-1) at (1,1) {$\vdots$};
        \node (0) at (1,0) {};
        \node (1l) at (2,-1) {};
        \node (1m) at (1,-1) {};
        \node (1r) at (0,-1) {};
        \node (2l) at (1,-2) {};
        \node (2m) at (0,-2) {};
        \node (2r) at (-1,-2) {};
        \node (3l) at (2,-3) {};
        \node (3m) at (1,-3) {};
        \node (3r) at (0,-3) {};
        \node (4m) at (1,-4) {$\vdots$};
        \draw (1l.center)--(0.center)--(1r.center);
        \draw (1m.center)--(0.center);
        \draw (2r.center)--(1r.center);
        \draw (2m.center)--(1r.center);
        \draw (2l.center)--(1r.center);
        \draw (3l.center)--(2l.center);
        \draw (3m.center)--(2l.center);
        \draw (3r.center)--(2l.center);
        \node (-1') at (-3,1) {$\vdots$};
        \node (0') at (-3,0) {};
        \node (1l') at (-2,-1) {};
        \node (1m') at (-3,-1) {};
        \node (1r') at (-4,-1) {};
        \node (2l') at (-3,-2) {};
        \node (2m') at (-4,-2) {};
        \node (2r') at (-5,-2) {};
        \node (2L') at (-3,-2) {};
        \node (2R') at (-5,-2) {};
        \node (3l') at (-2,-3) {};
        \node (3m') at (-3,-3) {};
        \node (3r') at (-4,-3) {};
        \node (4m') at (-3,-4) {$\vdots$};
        \draw (1l'.center)--(0'.center)--(1r'.center);
        \draw (1m'.center)--(0'.center);
        \draw (2r'.center)--(1r'.center);
        \draw (2m'.center)--(1r'.center);
        \draw (2l'.center)--(1r'.center);
        \draw (3l'.center)--(2l'.center);
        \draw (3m'.center)--(2l'.center);
        \draw (3r'.center)--(2l'.center);
    \end{tikzpicture}}=\left(\raisebox{-\height/2}{\begin{tikzpicture}[scale=.5]
        \node (-1) at (0,1) {$\vdots$};
        \node (0) at (0,0) {};
        \node (1l) at (-1,-1) {};
        \node (1m) at (0,-1) {};
        \node (1r) at (1,-1) {};
        \node (2l) at (-.25,-2) {};
        \node (2r) at (.25,-2) {};
        \node (2L) at (-1.25,-2) {};
        \node (2M) at (-1,-2) {};
        \node (2R) at (-.75,-2) {};
        \draw (1l.center)--(0.center)--(1r.center);
        \draw (1m.center)--(0.center);
        \draw (2r.center)--(1m.center);
        \draw (2l.center)--(1m.center);
        \draw (2L.center)--(1l.center);
        \draw (2M.center)--(1l.center);
        \draw (2R.center)--(1l.center);
        \node (-1') at (3,1) {$\vdots$};
        \node (0') at (3,0) {};
        \node (1l') at (2,-1) {};
        \node (1m') at (3,-1) {};
        \node (1r') at (4,-1) {};
        \node (2l') at (2.75,-2) {};
        \node (2m') at (3,-2) {};
        \node (2r') at (3.25,-2) {};
        \node (2L') at (1.75,-2) {};
        \node (2R') at (2.25,-2) {};
        \draw (1l'.center)--(0'.center)--(1r'.center);
        \draw (1m'.center)--(0'.center);
        \draw (2r'.center)--(1m'.center);
        \draw (2m'.center)--(1m'.center);
        \draw (2l'.center)--(1m'.center);
        \draw (2L'.center)--(1l'.center);
        \draw (2R'.center)--(1l'.center);
    \end{tikzpicture}}\right)\left(\raisebox{-\height/2}{\begin{tikzpicture}[scale=.5]
        \node (-1) at (0,1) {$\vdots$};
        \node (0) at (0,0) {};
        \node (1l) at (-1,-1) {};
        \node (1m) at (0,-1) {};
        \node (1r) at (1,-1) {};
        \node (2l) at (-1,-2) {};
        \node (2m) at (0,-2) {};
        \node (2r) at (1,-2) {};
        \node (3l) at (-2,-3) {};
        \node (3m) at (-1,-3) {};
        \node (3r) at (0,-3) {};
        \node (4m) at (-1,-4) {$\vdots$};
        \draw (1l.center)--(0.center)--(1r.center);
        \draw (1m.center)--(0.center);
        \draw (2r.center)--(1m.center);
        \draw (2m.center)--(1m.center);
        \draw (2l.center)--(1m.center);
        \draw (3l.center)--(2l.center);
        \draw (3m.center)--(2l.center);
        \draw (3r.center)--(2l.center);
        \node (-1') at (3,1) {$\vdots$};
        \node (0') at (3,0) {};
        \node (1l') at (2,-1) {};
        \node (1m') at (3,-1) {};
        \node (1r') at (4,-1) {};
        \node (2l') at (2,-2) {};
        \node (2m') at (3,-2) {};
        \node (2r') at (4,-2) {};
        \node (3l') at (1,-3) {};
        \node (3m') at (2,-3) {};
        \node (3r') at (3,-3) {};
        \node (4m') at (2,-4) {$\vdots$};
        \draw (1l'.center)--(0'.center)--(1r'.center);
        \draw (1m'.center)--(0'.center);
        \draw (2r'.center)--(1m'.center);
        \draw (2m'.center)--(1m'.center);
        \draw (2l'.center)--(1m'.center);
        \draw (3l'.center)--(2l'.center);
        \draw (3m'.center)--(2l'.center);
        \draw (3r'.center)--(2l'.center);
    \end{tikzpicture}}\right)\left(\raisebox{-\height/2}{\begin{tikzpicture}[scale=.5]
        \node (-1) at (0,1) {$\vdots$};
        \node (0) at (0,0) {};
        \node (1l) at (-1,-1) {};
        \node (1m) at (0,-1) {};
        \node (1r) at (1,-1) {};
        \node (2l) at (-.25,-2) {};
        \node (2m) at (0,-2) {};
        \node (2r) at (.25,-2) {};
        \node (2L) at (-1.25,-2) {};
        \node (2R) at (-.75,-2) {};
        \draw (1l.center)--(0.center)--(1r.center);
        \draw (1m.center)--(0.center);
        \draw (2r.center)--(1m.center);
        \draw (2m.center)--(1m.center);
        \draw (2l.center)--(1m.center);
        \draw (2L.center)--(1l.center);
        \draw (2R.center)--(1l.center);
        \node (-1') at (3,1) {$\vdots$};
        \node (0') at (3,0) {};
        \node (1l') at (2,-1) {};
        \node (1m') at (3,-1) {};
        \node (1r') at (4,-1) {};
        \node (2l') at (2.75,-2) {};
        \node (2r') at (3.25,-2) {};
        \node (2L') at (1.75,-2) {};
        \node (2M') at (2,-2) {};
        \node (2R') at (2.25,-2) {};
        \draw (1l'.center)--(0'.center)--(1r'.center);
        \draw (1m'.center)--(0'.center);
        \draw (2r'.center)--(1m'.center);
        \draw (2l'.center)--(1m'.center);
        \draw (2L'.center)--(1l'.center);
        \draw (2M'.center)--(1l'.center);
        \draw (2R'.center)--(1l'.center);
    \end{tikzpicture}}\right)$$
    $$\raisebox{-\height/2}{\begin{tikzpicture}[scale=.5]
        \node (-1) at (-1,1) {$\vdots$};
        \node (0) at (-1,0) {};
        \node (1l) at (-2,-1) {};
        \node (1m) at (-1,-1) {};
        \node (1r) at (0,-1) {};
        \node (2l) at (-1,-2) {};
        \node (2m) at (0,-2) {};
        \node (2r) at (1,-2) {};
        \node (3l) at (-2,-3) {};
        \node (3m) at (-1,-3) {};
        \node (3r) at (0,-3) {};
        \node (4m) at (-1,-4) {$\vdots$};
        \draw (1l.center)--(0.center)--(1r.center);
        \draw (1m.center)--(0.center);
        \draw (2r.center)--(1r.center);
        \draw (2m.center)--(1r.center);
        \draw (2l.center)--(1r.center);
        \draw (3l.center)--(2l.center);
        \draw (3m.center)--(2l.center);
        \draw (3r.center)--(2l.center);
        \node (-1') at (3,1) {$\vdots$};
        \node (0') at (3,0) {};
        \node (1l') at (2,-1) {};
        \node (1m') at (3,-1) {};
        \node (1r') at (4,-1) {};
        \node (2l') at (3,-2) {};
        \node (2m') at (4,-2) {};
        \node (2r') at (5,-2) {};
        \node (2L') at (3,-2) {};
        \node (2R') at (5,-2) {};
        \node (3l') at (2,-3) {};
        \node (3m') at (3,-3) {};
        \node (3r') at (4,-3) {};
        \node (4m') at (3,-4) {$\vdots$};
        \draw (1l'.center)--(0'.center)--(1r'.center);
        \draw (1m'.center)--(0'.center);
        \draw (2r'.center)--(1r'.center);
        \draw (2m'.center)--(1r'.center);
        \draw (2l'.center)--(1r'.center);
        \draw (3l'.center)--(2l'.center);
        \draw (3m'.center)--(2l'.center);
        \draw (3r'.center)--(2l'.center);
    \end{tikzpicture}}=\left(\raisebox{-\height/2}{\begin{tikzpicture}[scale=.5]
        \node (-1) at (0,1) {$\vdots$};
        \node (0) at (0,0) {};
        \node (1l) at (1,-1) {};
        \node (1m) at (0,-1) {};
        \node (1r) at (-1,-1) {};
        \node (2l) at (.25,-2) {};
        \node (2m) at (0,-2) {};
        \node (2r) at (-.25,-2) {};
        \node (2L) at (1.25,-2) {};
        \node (2R) at (.75,-2) {};
        \draw (1l.center)--(0.center)--(1r.center);
        \draw (1m.center)--(0.center);
        \draw (2r.center)--(1m.center);
        \draw (2m.center)--(1m.center);
        \draw (2l.center)--(1m.center);
        \draw (2L.center)--(1l.center);
        \draw (2R.center)--(1l.center);
        \node (-1') at (-3,1) {$\vdots$};
        \node (0') at (-3,0) {};
        \node (1l') at (-2,-1) {};
        \node (1m') at (-3,-1) {};
        \node (1r') at (-4,-1) {};
        \node (2l') at (-2.75,-2) {};
        \node (2r') at (-3.25,-2) {};
        \node (2L') at (-1.75,-2) {};
        \node (2M') at (-2,-2) {};
        \node (2R') at (-2.25,-2) {};
        \draw (1l'.center)--(0'.center)--(1r'.center);
        \draw (1m'.center)--(0'.center);
        \draw (2r'.center)--(1m'.center);
        \draw (2l'.center)--(1m'.center);
        \draw (2L'.center)--(1l'.center);
        \draw (2M'.center)--(1l'.center);
        \draw (2R'.center)--(1l'.center);
    \end{tikzpicture}}\right)\left(\raisebox{-\height/2}{\begin{tikzpicture}[scale=.5]
        \node (-1) at (0,1) {$\vdots$};
        \node (0) at (0,0) {};
        \node (1l) at (1,-1) {};
        \node (1m) at (0,-1) {};
        \node (1r) at (-1,-1) {};
        \node (2l) at (1,-2) {};
        \node (2m) at (0,-2) {};
        \node (2r) at (-1,-2) {};
        \node (3l) at (2,-3) {};
        \node (3m) at (1,-3) {};
        \node (3r) at (0,-3) {};
        \node (4m) at (1,-4) {$\vdots$};
        \draw (1l.center)--(0.center)--(1r.center);
        \draw (1m.center)--(0.center);
        \draw (2r.center)--(1m.center);
        \draw (2m.center)--(1m.center);
        \draw (2l.center)--(1m.center);
        \draw (3l.center)--(2l.center);
        \draw (3m.center)--(2l.center);
        \draw (3r.center)--(2l.center);
        \node (-1') at (-3,1) {$\vdots$};
        \node (0') at (-3,0) {};
        \node (1l') at (-2,-1) {};
        \node (1m') at (-3,-1) {};
        \node (1r') at (-4,-1) {};
        \node (2l') at (-2,-2) {};
        \node (2m') at (-3,-2) {};
        \node (2r') at (-4,-2) {};
        \node (3l') at (-1,-3) {};
        \node (3m') at (-2,-3) {};
        \node (3r') at (-3,-3) {};
        \node (4m') at (-2,-4) {$\vdots$};
        \draw (1l'.center)--(0'.center)--(1r'.center);
        \draw (1m'.center)--(0'.center);
        \draw (2r'.center)--(1m'.center);
        \draw (2m'.center)--(1m'.center);
        \draw (2l'.center)--(1m'.center);
        \draw (3l'.center)--(2l'.center);
        \draw (3m'.center)--(2l'.center);
        \draw (3r'.center)--(2l'.center);
    \end{tikzpicture}}\right)\left(\raisebox{-\height/2}{\begin{tikzpicture}[scale=.5]
        \node (-1) at (0,1) {$\vdots$};
        \node (0) at (0,0) {};
        \node (1l) at (1,-1) {};
        \node (1m) at (0,-1) {};
        \node (1r) at (-1,-1) {};
        \node (2l) at (.25,-2) {};
        \node (2r) at (-.25,-2) {};
        \node (2L) at (1.25,-2) {};
        \node (2M) at (1,-2) {};
        \node (2R) at (.75,-2) {};
        \draw (1l.center)--(0.center)--(1r.center);
        \draw (1m.center)--(0.center);
        \draw (2r.center)--(1m.center);
        \draw (2l.center)--(1m.center);
        \draw (2L.center)--(1l.center);
        \draw (2M.center)--(1l.center);
        \draw (2R.center)--(1l.center);
        \node (-1') at (-3,1) {$\vdots$};
        \node (0') at (-3,0) {};
        \node (1l') at (-2,-1) {};
        \node (1m') at (-3,-1) {};
        \node (1r') at (-4,-1) {};
        \node (2l') at (-2.75,-2) {};
        \node (2m') at (-3,-2) {};
        \node (2r') at (-3.25,-2) {};
        \node (2L') at (-1.75,-2) {};
        \node (2R') at (-2.25,-2) {};
        \draw (1l'.center)--(0'.center)--(1r'.center);
        \draw (1m'.center)--(0'.center);
        \draw (2r'.center)--(1m'.center);
        \draw (2m'.center)--(1m'.center);
        \draw (2l'.center)--(1m'.center);
        \draw (2L'.center)--(1l'.center);
        \draw (2R'.center)--(1l'.center);
    \end{tikzpicture}}\right)$$
    \caption{More conjugations to save the day.}
    \label{mesconjs}
    \end{figure}
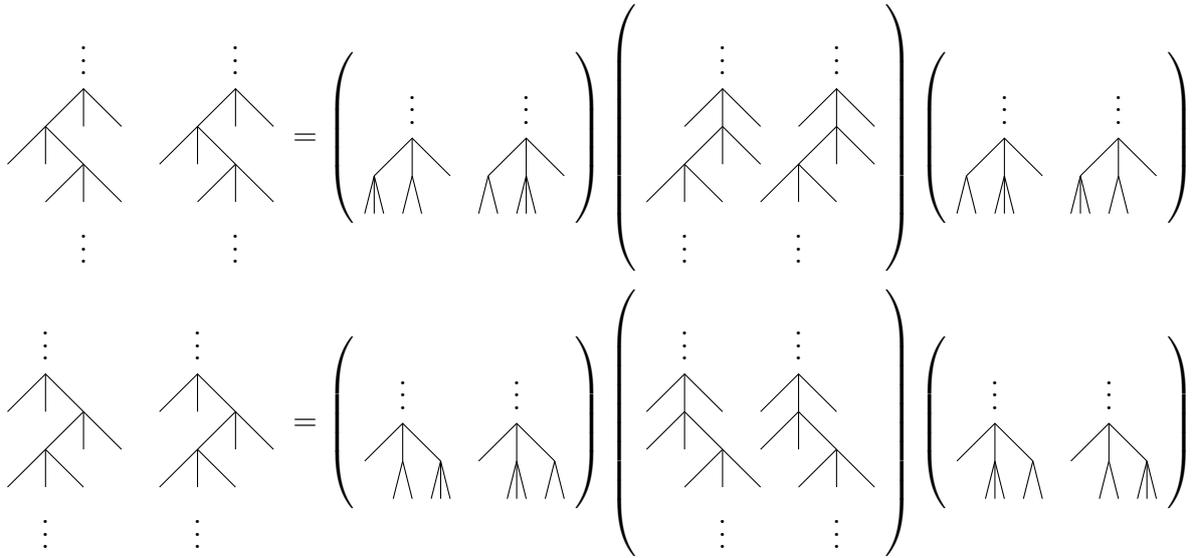

    The only centipedes left are those whose carets always hang from the left or always hang from the right, all the way to the hip. With $n$ torso pieces, there are only four of them (and their inverses), two of which are in the family we want to prove generates. The other two can be written in terms of stick bugs as seen in \Cref{quasistickbugs}.
    \begin{figure}[H]
    $$\raisebox{-\height/2}{\begin{tikzpicture}[scale=.6]
        \node (0) at (0,0) {};
        \node (1l) at (-1,-1) {};
        \node (1m) at (0,-1) {};
        \node (1r) at (1,-1) {};
        \node (2l) at (-2,-2) {};
        \node (2m) at (-1,-2) {};
        \node (2r) at (0,-2) {};
        \node (nl) at (-3,-3) {};
        \node (ml) at (-4,-4) {};
        \node (mm) at (-3,-4) {};
        \node (mr) at (-2,-4) {};
        \node (fl) at (-3.25,-5) {};
        \node (fm) at (-3,-5) {};
        \node (fr) at (-2.75,-5) {};
        \node (Fl) at (-2.25,-5) {};
        \node (Fr) at (-1.75,-5) {};
        \node (.) at (-2.5,-2.3) {$\iddots$};
        \draw (2l.center)--(1l.center)--(0.center)--(1r.center);
        \draw (1m.center)--(0.center);
        \draw (2r.center)--(1l.center);
        \draw (2m.center)--(1l.center);
        \draw (ml.center)--(nl.center);
        \draw (mm.center)--(nl.center);
        \draw (mr.center)--(nl.center);
        \draw (fl.center)--(mm.center);
        \draw (fm.center)--(mm.center);
        \draw (fr.center)--(mm.center);
        \draw (Fl.center)--(mr.center);
        \draw (Fr.center)--(mr.center);
        \draw[decoration={brace,mirror,raise=5pt},decorate] (0.center)--node[above left=4pt and 4pt] {$n$} (nl.center);
        \node (0') at (4,0) {};
        \node (1l') at (3,-1) {};
        \node (1m') at (4,-1) {};
        \node (1r') at (5,-1) {};
        \node (2l') at (2,-2) {};
        \node (2m') at (3,-2) {};
        \node (2r') at (4,-2) {};
        \node (nl') at (1,-3) {};
        \node (ml') at (0,-4) {};
        \node (mm') at (1,-4) {};
        \node (mr') at (2,-4) {};
        \node (fl') at (.75,-5) {};
        \node (fr') at (1.25,-5) {};
        \node (Fl') at (1.75,-5) {};
        \node (Fm') at (2,-5) {};
        \node (Fr') at (2.25,-5) {};
        \node (.') at (1.5,-2.3) {$\iddots$};
        \draw (2l'.center)--(1l'.center)--(0'.center)--(1r'.center);
        \draw (1m'.center)--(0'.center);
        \draw (2r'.center)--(1l'.center);
        \draw (2m'.center)--(1l'.center);
        \draw (ml'.center)--(nl'.center);
        \draw (mm'.center)--(nl'.center);
        \draw (mr'.center)--(nl'.center);
        \draw (fl'.center)--(mm'.center);
        \draw (fr'.center)--(mm'.center);
        \draw (Fl'.center)--(mr'.center);
        \draw (Fm'.center)--(mr'.center);
        \draw (Fr'.center)--(mr'.center);
        \draw[decoration={brace,mirror,raise=5pt},decorate] (0'.center)--node[above left=4pt and 4pt] {$n$} (nl'.center);
    \end{tikzpicture}}=l_n^{-1}l_{n-1}^{-1}l_nl_{n-1}=l_n^{-1}l^{-1}l_nl$$
    $$\raisebox{-\height/2}{\begin{tikzpicture}[scale=.6]
        \node (0) at (0,0) {};
        \node (1l) at (1,-1) {};
        \node (1m) at (0,-1) {};
        \node (1r) at (-1,-1) {};
        \node (2l) at (2,-2) {};
        \node (2m) at (1,-2) {};
        \node (2r) at (0,-2) {};
        \node (nl) at (3,-3) {};
        \node (ml) at (4,-4) {};
        \node (mm) at (3,-4) {};
        \node (mr) at (2,-4) {};
        \node (fl) at (3.25,-5) {};
        \node (fm) at (3,-5) {};
        \node (fr) at (2.75,-5) {};
        \node (Fl) at (2.25,-5) {};
        \node (Fr) at (1.75,-5) {};
        \node (.) at (2.5,-2.3) {$\ddots$};
        \draw (2l.center)--(1l.center)--(0.center)--(1r.center);
        \draw (1m.center)--(0.center);
        \draw (2r.center)--(1l.center);
        \draw (2m.center)--(1l.center);
        \draw (ml.center)--(nl.center);
        \draw (mm.center)--(nl.center);
        \draw (mr.center)--(nl.center);
        \draw (fl.center)--(mm.center);
        \draw (fm.center)--(mm.center);
        \draw (fr.center)--(mm.center);
        \draw (Fl.center)--(mr.center);
        \draw (Fr.center)--(mr.center);
        \draw[decoration={brace,raise=5pt},decorate] (0.center)--node[above right=4pt and 4pt] {$n$} (nl.center);
        \node (0') at (-4,0) {};
        \node (1l') at (-3,-1) {};
        \node (1m') at (-4,-1) {};
        \node (1r') at (-5,-1) {};
        \node (2l') at (-2,-2) {};
        \node (2m') at (-3,-2) {};
        \node (2r') at (-4,-2) {};
        \node (nl') at (-1,-3) {};
        \node (ml') at (0,-4) {};
        \node (mm') at (-1,-4) {};
        \node (mr') at (-2,-4) {};
        \node (fl') at (-.75,-5) {};
        \node (fr') at (-1.25,-5) {};
        \node (Fl') at (-1.75,-5) {};
        \node (Fm') at (-2,-5) {};
        \node (Fr') at (-2.25,-5) {};
        \node (.') at (-1.5,-2.3) {$\ddots$};
        \draw (2l'.center)--(1l'.center)--(0'.center)--(1r'.center);
        \draw (1m'.center)--(0'.center);
        \draw (2r'.center)--(1l'.center);
        \draw (2m'.center)--(1l'.center);
        \draw (ml'.center)--(nl'.center);
        \draw (mm'.center)--(nl'.center);
        \draw (mr'.center)--(nl'.center);
        \draw (fl'.center)--(mm'.center);
        \draw (fr'.center)--(mm'.center);
        \draw (Fl'.center)--(mr'.center);
        \draw (Fm'.center)--(mr'.center);
        \draw (Fr'.center)--(mr'.center);
        \draw[decoration={brace,raise=5pt},decorate] (0'.center)--node[above right=4pt and 4pt] {$n$} (nl'.center);
    \end{tikzpicture}}=r_{n-1}r_nr_{n-1}^{-1}r_n^{-1}=rr_nr^{-1}r_n^{-1}$$
    \caption{The last remaining centipedes written as a product of stickbugs (in two different ways).}
    \label{quasistickbugs}
    \end{figure}
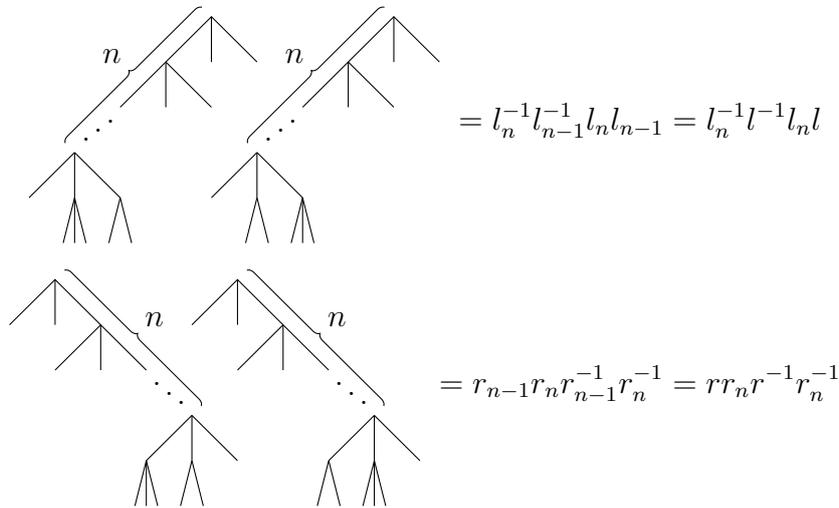

    This proves that the stick bugs generate the centipedes and hence $F\left(\frac32\right)$.
    \end{prf}

    \section{Finite generation}\label{finitegen}
    We are close to proving \Cref{fg}. Namely, we need to prove:

    \begin{thm}\label{lrgen}
        $l$ and $r$ generate the stick bugs.
    \end{thm}

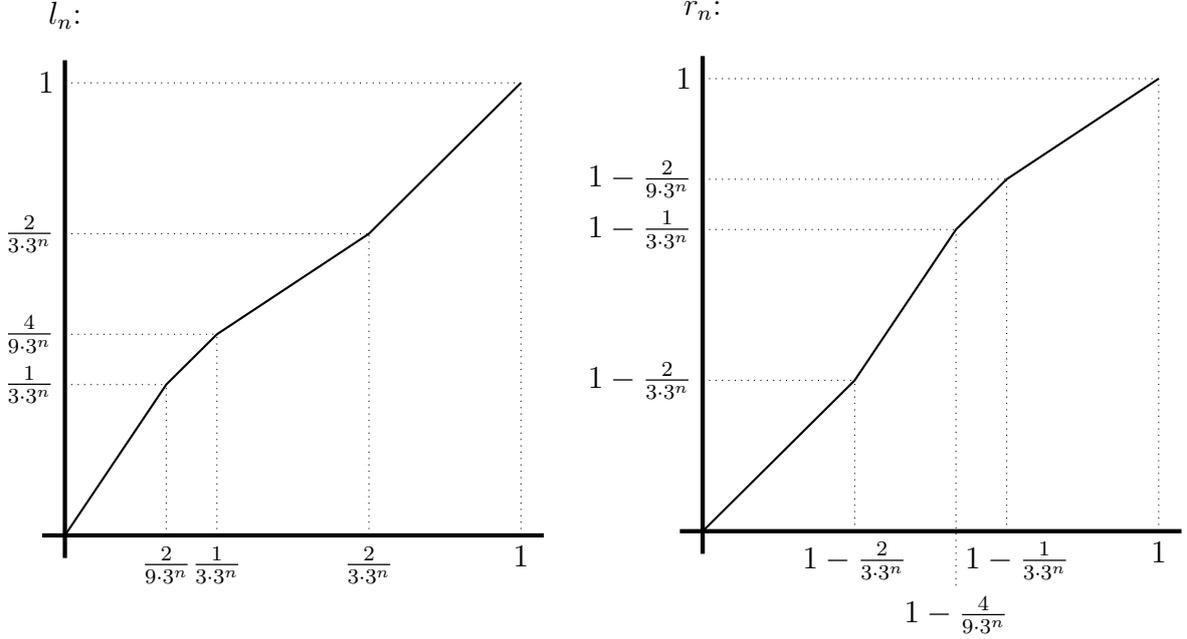
\begin{figure}[H]
        $$\raisebox{-\height}{\begin{tikzpicture}[scale=6]
            \node (l) at (0,1.15) {$l_n\text{:}$};
            \draw[thick] (0,0) -- (2/9,1/3) -- (1/3,4/9) -- (2/3,2/3) -- (1,1);
            \draw[ultra thick] (-.05,0) -- (1.05,0);
            \draw[ultra thick] (0,-.05) -- (0,1.05);
            \draw[dotted] (2/9,0) -- (2/9,1/3);
            \draw[dotted] (1/3,0) -- (1/3,4/9);
            \draw[dotted] (2/3,0) -- (2/3,2/3);
            \draw[dotted] (1,0) -- (1,1);
            \node[anchor=north] (x2/9) at (2/9,0) {$\frac2{9\cdot3^n}$};
            \node[anchor=north] (x1/3) at (1/3,0) {$\frac1{3\cdot3^n}$};
            \node[anchor=north] (x2/3) at (2/3,0) {$\frac2{3\cdot3^n}$};
            \node[anchor=north] (x1) at (1,0) {$1$};
            \draw[dotted] (0,1/3) -- (2/9,1/3);
            \draw[dotted] (0,4/9) -- (1/3,4/9);
            \draw[dotted] (0,2/3) -- (2/3,2/3);
            \draw[dotted] (0,1) -- (1,1);
            \node[anchor=east] (y1/3) at (0,1/3) {$\frac1{3\cdot3^n}$};
            \node[anchor=east] (y4/9) at (0,4/9) {$\frac4{9\cdot3^n}$};
            \node[anchor=east] (y2/3) at (0,2/3) {$\frac2{3\cdot3^n}$};
            \node[anchor=east] (y1) at (0,1) {$1$};
        \end{tikzpicture}}\quad\raisebox{-\height}{\begin{tikzpicture}[scale=6]
            \node (r) at (0,1.15) {$r_n\text{:}$};
            \draw[thick] (0,0) -- (1/3,1/3) -- (5/9,2/3) -- (2/3,7/9) -- (1,1);
            \draw[ultra thick] (-.05,0) -- (1.05,0);
            \draw[ultra thick] (0,-.05) -- (0,1.05);
            \draw[dotted] (1-1/3,0) -- (1-1/3,1-2/9);
            \draw[dotted] (1-4/9,-.12) -- (1-4/9,1-1/3);
            \draw[dotted] (1-2/3,0) -- (1-2/3,1-2/3);
            \draw[dotted] (1,0) -- (1,1);
            \node[anchor=north] (x1-1/3) at (1-1/3,0) {$\ \ 1-\frac1{3\cdot3^n}$};
            \node[anchor=north] (x1-4/9) at (1-4/9,-.12) {$1-\frac4{9\cdot3^n}$};
            \node[anchor=north] (x1-2/3) at (1-2/3,0) {$1-\frac2{3\cdot3^n}$};
            \node[anchor=north] (x1) at (1,0) {$1$};
            \draw[dotted] (0,1-2/9) -- (1-1/3,1-2/9);
            \draw[dotted] (0,1-1/3) -- (1-4/9,1-1/3);
            \draw[dotted] (0,1-2/3) -- (1-2/3,1-2/3);
            \draw[dotted] (0,1) -- (1,1);
            \node[anchor=east] (y1-2/9) at (0,1-2/9) {$1-\frac2{9\cdot3^n}$};
            \node[anchor=east] (y1-1/3) at (0,1-1/3) {$1-\frac1{3\cdot3^n}$};
            \node[anchor=east] (y1-2/3) at (0,1-2/3) {$1-\frac2{3\cdot3^n}$};
            \node[anchor=east] (y1) at (0,1) {$1$};
        \end{tikzpicture}}$$
        \caption{The graph of a stick bug.}
        \label{graficastickbug}
    \end{figure}

    In order to do this last bit, we will change tech and work with the functions for a little bit, instead of their tree-pair diagram. We show the graph of a general stick bug in \Cref{graficastickbug}, which can easily be derived from its tree-pair diagram.

    \begin{dfn}
    The \textit{support} $\supp(f)$ of an element $f$ in a Thompson-like group is the (topological) closure of the set of points $t$ such that $f(t)\neq t$.
    \end{dfn}

    It follows from the definition that not being in the support is a property closed by products and inverses. Therefore, if $w\in\left\langle f,g\right\rangle$, then $\supp(w)\subseteq\supp(f)\cup\supp(g)$.\\

    The proof that $l$ and $r$ generate the stick bugs will be done in two steps:

    \begin{prp}\label{prp1} For each $n\geq1$, there exist elements $X$ and $Y$ belonging to $\left\langle l,r\right\rangle$ such that $X\cdot l_n\cdot Y$ has support in $\left[\frac13,\frac23\right]$, and similarly for $r_n$.
    \end{prp}

    \begin{prp}\label{prp2} Every element with support in $\left[\frac13,\frac23\right]$ belongs to the span of $l$ and $r$.
    \end{prp}

    Before we prove \Cref{prp1}, we will need the following lemma:

    \begin{lmm}\label{lema} If $f$ is a function with support in $\left[\epsilon,\frac23\right]$, then $l^{-1}\cdot f\cdot l$ has support in $\left[l(\epsilon),\frac23\right]$.
    \end{lmm}
    \begin{prf}
        The support cannot exceed $\frac23$, and if $x<l(\epsilon)$, then $l^{-1}(x)<\epsilon$ follows from the monotonicity of $l^{-1}$, so $f(l^{-1}(x))=l^{-1}(x)$ and $l(f(l^{-1}(x)))=x$.
    \end{prf}

    \begin{prf}[ of \Cref{prp1}] We will prove that $l^{-k}\cdot(l^{-1}\cdot l_n)\cdot l^k$ has support in $\left[\frac13,\frac23\right]$ for some $k$, from which $X=l^{-k-1}$, $Y=l^k$. We have $\supp(l^{-1}\cdot l_n)\subseteq\supp(l)\cup\supp(l_n)=\left[0,\frac23\right]$, but since it is the identity near $0$, we have $\supp(l^{-1}\cdot l_n)\subseteq\left[\delta,\frac23\right]$ for some $\delta$ between $0$ and $\frac23$.\\

    By repeatedly applying \Cref{lema}, we can prove that $\supp(l^{-k}\cdot(l^{-1}\cdot l_n)\cdot l^k)\subseteq\left[\smash{\overbrace{l(l(\cdots l}^k(\delta)\cdots))},\frac23\right]$. Independently of the value $\delta>0$, the sequence $l(l(\cdots l(\delta)\cdots))$ is increasing and converges towards the fixed point $\frac23$ of $l$. So, it eventually exceeds $\frac13$.

    In a similar manner, we see that $\supp(r^k\cdot(r_n\cdot r^{-1})\cdot r^{-k})\subseteq\left[\frac13,\frac23\right]$, so the proposition is proven.
    \end{prf}

    Now, to prove \Cref{prp2}, we will need the following observation:

    \begin{obs}
    The elements with support inside $\left[\frac13,\frac23\right]$ form a subgroup of $F\left(\frac32\right)$ isomorphic to it. The isomorphism can be easily written in terms of tree-pair diagrams: if $(A,B)$ is a tree-pair diagram for $f\in F\left(\frac32\right)$ and $g$ is the image of $f$ under this isomorphism, then $(A',B')$ is a tree-pair diagram for $g$, where $A'$ and $B'$ are obtained by hanging $A$ and $B$, respectively, from the middle stem of a $3$-caret.
    \end{obs}

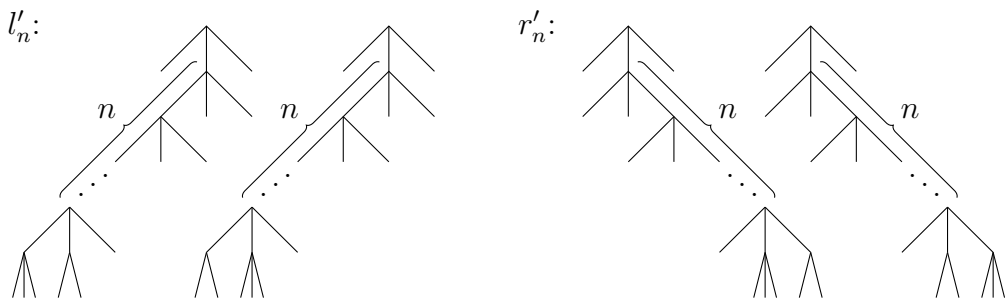
\begin{figure}[H]
    $$\raisebox{-\height/2}{\begin{tikzpicture}[scale=.6]
        \node (l) at (-4,1) {$l'_n\text{:}$};
        \node (-1) at (0,1) {};
        \node (0l) at (-1,0) {};
        \node (0) at (0,0) {};
        \node (0r) at (1,0) {};
        \node (1l) at (-1,-1) {};
        \node (1m) at (0,-1) {};
        \node (1r) at (1,-1) {};
        \node (2l) at (-2,-2) {};
        \node (2m) at (-1,-2) {};
        \node (2r) at (0,-2) {};
        \node (nl) at (-3,-3) {};
        \node (ml) at (-4,-4) {};
        \node (mm) at (-3,-4) {};
        \node (mr) at (-2,-4) {};
        \node (fl) at (-4.25,-5) {};
        \node (fm) at (-4,-5) {};
        \node (fr) at (-3.75,-5) {};
        \node (Fl) at (-3.25,-5) {};
        \node (Fr) at (-2.75,-5) {};
        \node (.) at (-2.5,-2.3) {$\iddots$};
        \draw (0l.center)--(-1.center)--(0r.center);
        \draw (-1.center)--(0.center);
        \draw (2l.center)--(1l.center)--(0.center)--(1r.center);
        \draw (1m.center)--(0.center);
        \draw (2r.center)--(1l.center);
        \draw (2m.center)--(1l.center);
        \draw (ml.center)--(nl.center);
        \draw (mm.center)--(nl.center);
        \draw (mr.center)--(nl.center);
        \draw (fl.center)--(ml.center);
        \draw (fm.center)--(ml.center);
        \draw (fr.center)--(ml.center);
        \draw (Fl.center)--(mm.center);
        \draw (Fr.center)--(mm.center);
        \draw[decoration={brace,mirror,raise=5pt},decorate] (0.center)--node[above left=4pt and 4pt] {$n$} (nl.center);
        \node (-1') at (4,1) {};
        \node (0l') at (3,0) {};
        \node (0') at (4,0) {};
        \node (0r') at (5,0) {};
        \node (1l') at (3,-1) {};
        \node (1m') at (4,-1) {};
        \node (1r') at (5,-1) {};
        \node (2l') at (2,-2) {};
        \node (2m') at (3,-2) {};
        \node (2r') at (4,-2) {};
        \node (nl') at (1,-3) {};
        \node (ml') at (0,-4) {};
        \node (mm') at (1,-4) {};
        \node (mr') at (2,-4) {};
        \node (fl') at (-.25,-5) {};
        \node (fr') at (.25,-5) {};
        \node (Fl') at (.75,-5) {};
        \node (Fm') at (1,-5) {};
        \node (Fr') at (1.25,-5) {};
        \node (.') at (1.5,-2.3) {$\iddots$};
        \draw (0l'.center)--(-1'.center)--(0r'.center);
        \draw (-1'.center)--(0'.center);
        \draw (2l'.center)--(1l'.center)--(0'.center)--(1r'.center);
        \draw (1m'.center)--(0'.center);
        \draw (2r'.center)--(1l'.center);
        \draw (2m'.center)--(1l'.center);
        \draw (ml'.center)--(nl'.center);
        \draw (mm'.center)--(nl'.center);
        \draw (mr'.center)--(nl'.center);
        \draw (fl'.center)--(ml'.center);
        \draw (fr'.center)--(ml'.center);
        \draw (Fl'.center)--(mm'.center);
        \draw (Fm'.center)--(mm'.center);
        \draw (Fr'.center)--(mm'.center);
        \draw[decoration={brace,mirror,raise=5pt},decorate] (0'.center)--node[above left=4pt and 4pt] {$n$} (nl'.center);
    \end{tikzpicture}}\qquad\raisebox{-\height/2}{\begin{tikzpicture}[scale=.6]
        \node (l) at (-6,1) {$r'_n\text{:}$};
        \node (-1) at (0,1) {};
        \node (0l) at (-1,0) {};
        \node (0) at (0,0) {};
        \node (0r) at (1,0) {};
        \node (1l) at (1,-1) {};
        \node (1m) at (0,-1) {};
        \node (1r) at (-1,-1) {};
        \node (2l) at (2,-2) {};
        \node (2m) at (1,-2) {};
        \node (2r) at (0,-2) {};
        \node (nl) at (3,-3) {};
        \node (ml) at (4,-4) {};
        \node (mm) at (3,-4) {};
        \node (mr) at (2,-4) {};
        \node (fl) at (4.25,-5) {};
        \node (fm) at (4,-5) {};
        \node (fr) at (3.75,-5) {};
        \node (Fl) at (3.25,-5) {};
        \node (Fr) at (2.75,-5) {};
        \node (.) at (2.5,-2.3) {$\ddots$};
        \draw (0l.center)--(-1.center)--(0r.center);
        \draw (0.center)--(-1.center);
        \draw (2l.center)--(1l.center)--(0.center)--(1r.center);
        \draw (1m.center)--(0.center);
        \draw (2r.center)--(1l.center);
        \draw (2m.center)--(1l.center);
        \draw (ml.center)--(nl.center);
        \draw (mm.center)--(nl.center);
        \draw (mr.center)--(nl.center);
        \draw (fl.center)--(ml.center);
        \draw (fm.center)--(ml.center);
        \draw (fr.center)--(ml.center);
        \draw (Fl.center)--(mm.center);
        \draw (Fr.center)--(mm.center);
        \draw[decoration={brace,raise=5pt},decorate] (0.center)--node[above right=4pt and 4pt] {$n$} (nl.center);
        \node (-1') at (-4,1) {};
        \node (0l') at (-5,0) {};
        \node (0') at (-4,0) {};
        \node (0r') at (-3,0) {};
        \node (1l') at (-3,-1) {};
        \node (1m') at (-4,-1) {};
        \node (1r') at (-5,-1) {};
        \node (2l') at (-2,-2) {};
        \node (2m') at (-3,-2) {};
        \node (2r') at (-4,-2) {};
        \node (nl') at (-1,-3) {};
        \node (ml') at (0,-4) {};
        \node (mm') at (-1,-4) {};
        \node (mr') at (-2,-4) {};
        \node (fl') at (.25,-5) {};
        \node (fr') at (-.25,-5) {};
        \node (Fl') at (-.75,-5) {};
        \node (Fm') at (-1,-5) {};
        \node (Fr') at (-1.25,-5) {};
        \node (.') at (-1.5,-2.3) {$\ddots$};
        \draw (0l'.center)--(-1'.center)--(0r'.center);
        \draw (0'.center)--(-1'.center);
        \draw (2l'.center)--(1l'.center)--(0'.center)--(1r'.center);
        \draw (1m'.center)--(0'.center);
        \draw (2r'.center)--(1l'.center);
        \draw (2m'.center)--(1l'.center);
        \draw (ml'.center)--(nl'.center);
        \draw (mm'.center)--(nl'.center);
        \draw (mr'.center)--(nl'.center);
        \draw (fl'.center)--(ml'.center);
        \draw (fr'.center)--(ml'.center);
        \draw (Fl'.center)--(mm'.center);
        \draw (Fm'.center)--(mm'.center);
        \draw (Fr'.center)--(mm'.center);
        \draw[decoration={brace,raise=5pt},decorate] (0'.center)--node[above right=4pt and 4pt] {$n$} (nl'.center);
    \end{tikzpicture}}$$
    \caption{The image of the stick bugs' tree-pair diagram under the isomorphism.}
    \label{pseudostickbugs}
    \end{figure}

    A consequence of this observation, together with the fact that the stick bugs are generators (\Cref{stickgen}), is that the elements $\{l'_0,l'_1,l'_2,\dots,r'_0,r'_1,r'_2,\dots\}$, defined as \Cref{pseudostickbugs}, generate all functions with support in $\left[\frac13,\frac23\right]$, since they are the images of the stick bugs.

    \begin{prf}[ of \Cref{prp2}] Because of the above, it is sufficient to write $l'_n$ and $r'_n$ in terms of $l$ and $r$. If we repeatedly apply the decompositions in \Cref{redreçantcentpeus}, we see that all the factors are tree-pair diagrams where every torso and hip caret hangs from the middle stem of the previous one. These belong to the span of $l$ and $r$ because of \Cref{xvnested}.
    \end{prf}

    We are now able to prove that $l$ and $r$ generate the stick bugs:

    \begin{prf}[ of \Cref{lrgen}] Follows from \Cref{prp1} and \Cref{prp2}, passing $X$ and $Y$ to the other side.\end{prf}

    This reduces the generating family to a set of two elements, concluding the proof of $\Cref{fg}$.

\bibliography{pepsrefs}
\bibliographystyle{plain}

\end{document}